\documentclass[12pt]{amsart}
\usepackage{amscd,amssymb}
\newcommand{\Aut}{\operatorname{Aut}}
\newcommand{\bb}{\mathbf{b}}
\newcommand{\cft}{CohFT}
\newcommand{\cfts}{CohFTs}
\newcommand{\C}{\mathcal{C}}	% The Universal Curve
\newcommand{\CP}[1]{\mathbb{CP}^{#1}}
\newcommand{\chit}{\widetilde{\chi}} % \chi with all u_i=0 except u_1
\newcommand{\EE}{\mathcal{E}} % Euler field
\newcommand{\End}[1]{\mathcal{E}nd({#1})}
\newcommand{\G}{\mathcal{G}} 	% Set of Stable Graphs
\newcommand{\GG}[1]{\nc[\G_{#1}]}
\newcommand{\HH}{\mathcal{H}}	% Operad of Graphs modulo relations

\newcommand{\jf}{\mathbf{j}}
\renewcommand{\ll}{\mathbf{l}}    
\newcommand{\LL}{\mathcal{L}}	% Line Bundle
\newcommand{\mf}{\mathbf{m}}    
\newcommand{\M}{\overline{\MM}}	% Compactified Moduli Space
\newcommand{\MM}{\mathcal{M}}	% Uncompactified Moduli Space
\newcommand{\pf}{\mathbf{p}}
\newcommand{\R}{\mathcal{R}}	% relations in the operadic sense
\newcommand{\SC}{\mathcal{S}}
\renewcommand{\sf}{\mathbf{s}}
\newcommand{\tf}{\mathbf{t}}
\newcommand{\tr}{{\rm tr}}	% the trace operation

\newcommand{\dd}{\partial} % partial w.r.t. t_a, x
\newcommand{\delf}{\boldsymbol{\delta}}
\newcommand{\epsf}{\boldsymbol{\varepsilon}}
\newcommand{\kaf}{\boldsymbol{\kappa}} % fat kappa
\newcommand{\kao}{\widehat{\kappa}} % kappa classes on M_{g,n+1}
\newcommand{\psio}{\widehat{\psi}} % psi classes on M_{g,n+1}
\newcommand{\tauf}{\boldsymbol{\tau}}
\newcommand{\zef}{\boldsymbol{0}} % far zero

\newcommand{\nc}{{\mathbb{C}}} 	% complex numbers
\newcommand{\nq}{{\mathbb{Q}}} 	% rational numbers
\newcommand{\nz}{{\mathbb{Z}}}	% integers

\newcommand{\lf}{\boldsymbol{[}}
\newcommand{\rf}{\boldsymbol{]}}
\newcommand{\la}{\langle}  % left angle bracket
\newcommand{\ra}{\rangle}  % right angle bracket

\def\lift#1#2{
  \dimen0 = \unitlength 
  \multiply\dimen0 by #1 \divide \dimen0 by 2
  \dimen1 = \dimen0 
  \multiply \dimen1 by 7 \divide \dimen1 by 10
  \raise\dimen1
     \hbox{\hskip 0.3cm ${\vbox to \dimen0{}}$ \enspace #2}}

\newtheorem{thm}{Theorem}[section]
\newtheorem{lm}[thm]{Lemma}
\newtheorem{prop}[thm]{Proposition}
\newtheorem{crl}[thm]{Corollary}

\theoremstyle{definition}

\newtheorem{df}[thm]{Definition}

\theoremstyle{remark}

\newtheorem{nota}{Notation} 
\newtheorem{rem}{Remark}  %[section] 
\newtheorem{ack}{Acknowledgment}

\begin{document}

\title[Intersection Numbers]
{Intersection Numbers and Rank One Cohomological Field Theories in Genus One}

\author
[A. Kabanov]{Alexandre Kabanov}
\address
{Max Planck Institut f\"ur Mathematik, Gottfried Claren Str. 26, 53225
Bonn, GERMANY and Department of Mathematics, Michigan State
University, Wells Hall,
East Lansing, MI 48824-1027, USA}
\email{kabanov@math.msu.edu}
%\thanks{Research of the second author was supported by... }

\author
[T. Kimura]{Takashi Kimura}
\address
{Max Planck Institut f\"ur Mathematik, Gottfried Claren Str. 26, 53225
Bonn, GERMANY and Department of Mathematics, Boston University, Boston, MA
02215, USA} 
\email{kimura@math.bu.edu}
%\thanks{Research of the first author was supported by...}

\date{June 18, 1997}

\begin{abstract} We obtain a simple, recursive presentation of the
tautological ($\kappa$, $\psi$, and $\lambda$) classes on the moduli
space of curves in genus $0$ and $1$ in terms of boundary strata
(graphs). We derive differential equations for the generating
functions for their intersection numbers which allow us to prove a
simple relationship between the genus zero and genus one potentials.
As an application, we describe the moduli space of normalized, restricted,
rank one cohomological field theories in genus one in coordinates
which are additive under taking tensor products. Our results simplify
and generalize those of Kaufmann, Manin, and Zagier.
\end{abstract}

\maketitle

%	the introduction
% Last Modified 06/13/97 9:00pm by Sasha
% intro.tex

Recently, there has been a great deal of interest in the topology of
the moduli space of curves. Much of this interest has been due to the
important role that these spaces (and their cousins, the moduli space
of stable maps) play in the theory of Gromov-Witten invariants and
quantum cohomology \cite{KM,W,RT} whose origins in the physical
literature are called a topological gravity \cite{W}. They furnish
nontrivial examples of cohomological field theories (\cfts) \cite{KM,
Ma}, in genus zero (and conjecturally for higher genera). Often, this
structure is enough to completely determine the Gromov-Witten
invariants themselves. The moduli spaces of curves are endowed with
tautological classes whose generating functions for their associated
intersection numbers obey a system of differential equations which
often possess remarkable properties \cite{W,Ko}. In this paper, we
apply a mixture of algebraic geometry and combinatorics to find a
simple presentation of these classes in genus $0$ and $1$ to obtain a
generalization of some equations due to Witten and Dijkgraaf
\cite{W,Di}. These generating functions parameterize the potentials
associated to the space of all normalized, restricted, rank one
cohomological field theories in genus one and endow this space with
coordinates which are additive with respect to tensor product in the
category of \cfts. This paper is motivated by the work of Kaufmann,
Manin, and Zagier \cite{KMZ}. 

The moduli space of genus $g$ curves with $n$ marked points,
$\MM_{g,n}\,:=\, \{\,[C;x_1,x_2,\ldots,x_n]\,\}$, is the moduli space
of configurations of $n$ marked points, on a smooth, complex curve
(Riemann surface) $C$ of genus $g$. We assume throughout that the
stability condition $2-2g - n < 0$ is satisfied. This moduli space has
a compactification $\M_{g,n}$ (due to Deligne-Knudsen-Mumford) which
is the moduli space of stable curves of genus $g$ with $n$ marked
points where the boundary divisor $\M_{g,n}-\MM_{g,n}$ is the locus of
degenerate curves. The space $\M_{g,n}$ is a stratified, complex
orbifold (stack) of complex dimension $3g-3+n$ where each stratum is
indexed by a decorated graph (stable graph) which denotes the type of
degeneration that curves in that stratum have. A stable graph
represents a cohomology class on $\M_{g,n}$ by taking the closure of
the corresponding stratum and applying Poincar\'e duality to the
associated (rational) homology class.

The space $\M_{g,n}$ is endowed with tautological cohomology classes
whose study was initiated by Mumford \cite{Mu}. Let $\LL_i$ be the
line bundle over $\M_{g,n}$ whose fiber over a point
$[C;x_1,x_2,\ldots,x_n]$ is $T^*_{x_i} C$. Then $\psi_{(g,n),i} = c_1
(\LL_i)$, the first Chern class of $\LL_i$. The classes
$\kappa_{(g,n),i}$ in $H^{2i}(\M_{g,n},\nq)^{S_n}$ are defined by
$\kappa_{(g,n),i}\,:=\,\pi_*(c_1(\omega_{g,n}(D))^{i+1})$ where
$\omega_{g,n}$ is the cotangent bundle to the fibers of the
universal curve $\pi\,:\,\M_{g,n+1} \,\to\, \M_{g,n}$, $D$ is the sum
of the images of the canonical sections, and $\pi_*$ is fiber
integration \cite{AC}. Integrals of products of these classes
(intersection numbers) are of great geometric interest and are the
main object of study in this paper. In particular, $\frac{1}{2\pi^2}\,
\kappa_{(g,n),1}$ is the class of the Weil-Petersson symplectic form
\cite{AC}. Zograf in \cite{Z} obtained a recursion formula for the
classical Weil-Petersson volumes in genus zero. The intersection
numbers of the $\kappa$ classes studied in \cite{KMZ} are called
higher Weil-Petersson volumes.

In the first part of this paper we study a generating function
$H(\tf;\sf) \in \nc[[\tf,\sf]]$ which incorporates all intersection
numbers of the $\psi$ and $\kappa$ classes. Here $\tf = (t_0, t_1,
\ldots)$ and $\sf= (s_1, s_2, \ldots)$ are formal variables. The
function $H_g(\tf;\sf)$ denotes the summand of $H(\tf;\sf)$
corresponding to genus $g$. This function has the property that
$H(\tf,\zef) = F(\tf)$, the generating function for the $\psi$
intersection numbers defined by Witten in \cite{W,W2}. On the other
hand, setting $t_i=0$ for $i\ge 1$ gives a generating function for the
intersection numbers of the $\kappa$ classes. This function is closely
related to that considered in \cite{KMZ}. 

We give a simple, recursive presentation of powers of the $\kappa$ and
$\psi$ classes in terms of boundary strata in genus zero and genus one
and derive a simple system of differential equations for $H$ in genus
zero and one which completely determine those intersection numbers.
Taking appropriate limits in genus zero, we obtain equations due to
Witten \cite{W} for the $\psi$ classes and equations for the $\kappa$
classes which are equivalent of those in \cite{KMZ} and much simpler.
(A differential equation for the classical Weil-Petersson volumes in
genus zero was first obtained in \cite{Mat}.) This simplification
arises because our presentation of the $\kappa$ classes in terms of
boundary strata is simpler. Furthermore, the genus one equation can be
solved to obtain the relation
\begin{equation}
\label{our:rel}
H_1 = \frac{1}{24}\log H_0''',
\end{equation}
where $'$ denotes the partial derivative with respect to $t_0$, for
all values of $\mathbf{s}$ and $\mathbf{t}$ generalizing the result of
Dijkgraaf and Witten saying that $F_1= \frac{1}{24} \log F_0'''$
\cite{Di}. We were informed that \eqref{our:rel} was known to Zograf
\cite{Ma2} in the special case where $t_j = 0$ for all $j\geq 1$ and
$s_i=0$ for all $i\geq 2$.

Witten \cite{W} conjectured (and Kontsevich \cite{Ko} proved) that
$F(\mathbf{s})$ was the logarithm of a tau function of the KdV
hierarchy after rescaling the variables. This tau function was
completely characterized \cite{Di,KS,Lo} by being annihilated by a
sequence of differential operators $L_n$ for $n\geq -1$ satisfying
$[L_m,L_n] = (m-n) L_{m+n}$, the relations of the Virasoro
algebra.\footnote{This result is all the more interesting because of
the recent conjecture in \cite{E} which predicts a Virasoro algebra
playing a similar role in the the case where the target manifold is
nontrivial.} The equations corresponding to $L_{-1}$ and $L_0$ were
proven by Witten in \cite{W,W2} and are essentially the so-called
puncture and dilaton equations. We derive the analog of the puncture
and dilaton equations for intersection numbers of $\kappa$ and $\psi$
classes. These equations are valid for all genera and do not use the
presentation of the tautological classes in terms of boundary strata.
Solving these equations provides another proof of \eqref{our:rel}. It
is not clear which one of these approaches will prove to be most
useful in higher genera.

In the second part of our paper, we apply the previous to describe the
moduli space of normalized, restricted, rank one cohomological field
theories (\cfts) in genus one generalizing the results of Kaufmann,
Manin, and Zagier in genus zero \cite{KMZ}.

A (complete) \cft\ of rank $r$ \cite{KM,Ma} is an $r$ dimensional
vector space with metric $(V,h)$ together with a collection of linear
maps $H_\bullet(\M_{g,n})\,\to\,T^n V$ which are equivariant under the
action of the permutation groups and which satisfy some compatibility
conditions arising from inclusion of strata on $\M_{g,n}$ In the
language of Getzler and Kapranov \cite{GK}, the maps form a morphism
of modular operads. Restricting to $g=0$, a \cft\ $(V,h)$ is
equivalent to endowing $(V,h)$ with the structure of a (formal)
Frobenius manifold \cite{Du,Hi,Ma}. The most spectacular examples of
such theories arise when $(V,h)$ is the cohomology ring of certain
smooth, projective varieties with its intersection pairing and the
morphisms come from the Gromov-Witten invariants associated to the
manifold thereby endowing the cohomology ring with a deformed cup
product giving it the structure of quantum cohomology. In many cases,
{\sl e.g.}\ $\CP{n}$ or Grassmannians, the structure of a \cft\ is
strong enough to completely determine, recursively, the number of
rational curves in the manifold counted with multiplicity (see
\cite{KM,FP}).

A \cft\ in genus zero can be described in terms of a certain
generating function (potential) associated to the structure morphisms
which must satisfy the WDVV (Witten-Dijkgraaf-Verlinde-Verlinde)
equations (see \cite{FP,KM}). These equations encode the relations
between the boundary strata in $\M_{0,n}$ due to Keel in genus zero.
Recently, Getzler \cite{G} derived equations which are the analogs of
WDVV in genus one by proving new relations which plays an analogous role
in genus one to those of Keel in genus zero. His equation allowed him
to predict the elliptic Gromov-Witten invariants of $\CP{2}$ and
$\CP{3}$ (see also \cite{CH}).

Kaufmann, Manin, and Zagier proved \cite{KMZ} that the moduli space of
normalized, rank one, \cfts\ in genus zero has coordinates
$\mathbf{s}$ such that the tensor product in the category of \cfts\ is
additive in these coordinates. We prove that the moduli space of
normalized, rank one, restricted, \cfts\ in genus one has similar
coordinates $(\mathbf{s},u)$. Here the variable $u$ arises because in
genus one we need to introduce another tautological cohomology class
(called $\lambda$) due to Mumford. Our proof involves a mixture of
techniques from Kaufmann, Manin, and Zagier \cite{KMZ}, Getzler
\cite{G}, and results from the first part of this paper.

In section $1$, we review the geometry of the moduli space of stable
curves, its stratification in terms of stable graphs, tautological
cohomology classes, their intersection numbers, and associated
generating functions. In section 2, we obtain a simple presentation
for these classes in terms of stable graphs in genus $0$ and $1$ and
derive differential equations satisfied by the generating functions
associated to their intersection numbers. In section 3, we derive the
analogues of the puncture and dilaton equations. In section 4, we
write closed form expressions for these intersection numbers. In
section 5, we use analytic properties of the generating function to
prove an asymptotic formula for the Weil-Petersson volumes of the
moduli space of genus 1 curves as the number of punctures becomes very
large. Finally in section 6, we describe the moduli space of rank 1,
restricted, \cfts\ in genus one.

\begin{sloppypar}
\begin{ack}
We are greatful to the Max Planck Institut f\"ur Mathematik for their
financial support and for providing a wonderfully stimulating
atmosphere. We would like to thank R. Dijkgraaf, E. Getzler, and Yu.
Manin for useful conversations. We are grateful to J. Stasheff for his
comments on an earlier version of this paper. We would also like to
thank K. Belabas for his \TeX nical assistance and for providing the
music.
\end{ack}
\end{sloppypar}

%	moduli space of curves, graphs, taut classes, gen functions
% Last Modified 05/30/97 4pm by TK
\section{Moduli Space of Curves}
\label{kmz}

\begin{nota}
In this paper we always consider cohomology with the rational
coefficients: $H^\bullet(X)$ stands for $H^\bullet(X;\nq)$. We denote
the set $\{ 1, \ldots, n \}$ by $[n]$. If $I$ is a finite set we
denote its cardinality by $|I|$.
\end{nota}

\subsection{Basic Definitions}

Let $\MM_{g,n}$ be the moduli space of smooth curves of genus $g$ with
$n$ marked points, where $2g-2+n >0$, i.e.\ $\MM_{g,n} = \{\, [\Sigma
; x_1,x_2,\ldots,x_n ]\, \}$ where $\Sigma$ is a genus $g$ Riemann
surface and $x_1,x_2,\ldots,x_n$ are distinct marked points on
$\Sigma$. Two such configurations are equivalent if they are related
by a biholomorphic map. The moduli space $\MM_{g,n}$ has a natural
compactification due to Deligne, Knudsen, and Mumford denoted by
$\M_{g,n} = \{\,[C; x_1,x_2,\ldots, x_n] \,\}$ which is the moduli
space of stable curves of genus $g$ with $n$ punctures in which
$\MM_{g,n}$ sits as a dense open subset. The spaces $\M_{g,n}$ are
connected, compact, complex orbifolds (in fact, stacks) with complex
dimensions $3g-3+n$. The complement of $\MM_{g,n}$ in $\M_{g,n}$ is a
divisor with normal crossings and consists of those stable curves
which have double points.

The moduli space $\M_{g,n}$ forms the base of a universal family. Let
$\pi\,:\,\C_{g,n}\,\to\, \M_{g,n}$ be the \emph{universal curve} which
can be identified with $\C_{g,n} = \M_{g,n+1}$ where $\pi$ is the
projection obtained by forgetting the $(n+1)^{\text{st}}$ puncture and
followed by collapsing any resulting unstable irreducible components
of the curve, if any, to a point. The universal curve
$\C_{g,n}\,\to\,\M_{g,n}$ is furthermore endowed with canonical
sections $\sigma_1,\sigma_2,\ldots,\sigma_n$ such that $\sigma_i$ maps
$[C;x_1,x_2,\ldots,x_n]\,\mapsto\,[C';x_1',x_2',\ldots,x_{n+1}']$
where $C'$ is obtained from $C$ by attaching a three punctured sphere
to $x_i$ at one of its punctures to create a double point, then labeling
the remaining two punctures on the sphere $x_i'$ and $x_{n+1}'$, and
finally setting all other $x_j' = x_j$. The sections $\sigma_i$ are
well-defined since $\M_{0,3}$ is a point. The image of
$\sigma_i:\M_{g,n}\,\to\,\C_{g,n}$ gives rise to a divisor $D_{i}$ in
$\C_{g,n}$ for all $i=1,2,\ldots,n$.

\subsection{Natural Stratification} 

In the sequel it will be convenient to consider markings by arbitrary
finite sets rather than by just $[n]$. If $I$ is a finite set we
denote by $\M_{g,I} \cong \M_{g,|I|}$ the corresponding moduli space.

The natural stratification of $\M_{g,n}$ is best described in terms of
graphs, and therefore we start with fixing the notation concerning
graphs. We will consider only connected graphs. Each graph $\Gamma$ can be
described in terms of its set of vertices $V(\Gamma)$, set of edges
$E(\Gamma)$, and set of tails $S(\Gamma)$. Each edge has two endpoints
belonging to $V(\Gamma)$ which are allowed to be the same. Each tail
has only one endpoint. If $v\in V(\Gamma)$,we denote by $n(v)$ the
number of half-edges emanating from $v$, where each edge gives rise to
two half-edges, and each tail to one half-edge.

The natural stratification of $\M_{g,n}$ is determined by the type of
the degeneration of the curve representing a point in the moduli
space, and its strata can be labeled by stable graphs. A \emph{stable
graph} consists of a triple $(\Gamma,g,\mu)$, where $\Gamma$ is a
connected graph as above, $g: V(\Gamma) \to \nz_{\ge 0}$, and $\mu$ is
a bijection between $S(\Gamma)$ and a given set $I$. Moreover, one
requires that for each vertex $v$, the stability condition $2g(v)-2+
n(v)>0$ is satisfied. If $[C; x_1,x_2,\ldots, x_n]$ is a stable,
$n$-pointed curve one obtains the corresponding stable graph, called
its \emph{dual graph}, by collapsing each irreducible component to a
point (vertex), connecting any two vertices if their corresponding
components share a double point and attaching a tail to a vertex for
each marked point on that component.

We define the \emph{genus $g(\Gamma)$ of $\Gamma$} to be $b_1(\Gamma)+
\sum_{v\in V(\Gamma)} g(v)$, where $b_1(\Gamma)$ is the first Betti
number of $\Gamma$. We denote by $\G_{g,n}$ the set of the equivalence
classes of stable graphs of genus $g$ with $n$ tails labeled by $[n]$.
There is a natural action of the symmetric group $S_n$ on $\G_{g,n}$.
Associating a stable curve to its dual graph provides an
$S_n$-equivariant bijection between the strata of the natural
stratification of $\M_{g,n}$ and the elements of $\G_{g,n}$.

Let $\M_\Gamma$ be (the closure of) the moduli space of stable curves whose
dual graph is $\Gamma$. It is a closed irreducible subvariety of
codimension $|E(\Gamma)|$ of $\M_{g(\Gamma), S(\Gamma)}$. Moreover, it
is isomorphic to a quotient of the cartesian product
\[
\prod_{v\in V(\Gamma)} \M_{g(v),n(v)}
\]
by $\Aut(\Gamma)$, where the automorphisms of a stable graph
$(\Gamma,g,\mu)$ are required to preserve $g$ and $\mu$. This quotient
morphism can be made canonical if one creates a pair of labels for
each edge of $\Gamma$ and labels the $n(v)$ half-edges emanating from
$v$ by the corresponding elements of $S(\Gamma)$ with the labels
corresponding to the edges.

Each $\M_\Gamma$ determines the fundamental class, in the sense of
orbifolds, lying in $H^\bullet(\M_{g(\Gamma), n(\Gamma)})$. The pull
back of this class under the morphism $\pi: \M_{g,n+1} \to \M_{g,n}$
is represented by a subvariety of $\M_{g,n+1}$ corresponding to the
$|V(\Gamma)|$ graphs each of which is obtained by attaching a tail
numbered $n+1$ to a vertex of $\Gamma$. It is also easy to push down
these fundamental classes. Let $\M_{\Gamma'}$ represents an element of
$H^\bullet(\M_{g,n+1})$. The image of this element under $\pi_* :
H^{\bullet +2}(\C_{g,n})\,\to\,H^\bullet (\M_{g,n})$, induced by the
fiber integration, is zero if after removing the $(n+1)^{\text{st}}$
tail from $\Gamma'$ the graph remains stable. In the other case, when
the removal of the $(n+1)^{\text{st}}$ tail destabilizes $\Gamma'$,
the image is obtained by stabilization, i.e., contracting the edge
connecting the unstable vertex with the rest of the graph.

\subsection{Tautological Classes}
\label{tautological}

We will now describe three types of tautological cohomology classes
($\psi$, $\kappa$, and $\lambda$) associated to the universal curve.
Consider the universal curve $\C_{g,n}\,\to\,\M_{g,n}$. The cotangent
bundle to its fibers (in the orbifold sense) forms the holomorphic
line bundle $\omega_{g,n}$. Let $\LL_{(g,n),i}\,\to\,\M_{g,n}$ be
given by the pullback $\LL_{(g,n),i} = \sigma_i^* \omega_{(g,n)}$. The
tautological classes $\psi_{(g,n),i}$ in $H^2(\M_{g,n})$ are defined
by
\[
\psi_{(g,n),i} := c_1(\LL_{(g,n),i})
\]
where $c_1$ denotes the first Chern class.\footnote{The Chern classes
are in the sense of orbifolds and are therefore rational.}

The tautological classes $\kappa_{(g,n),i}$ in $H^{2i}(\M_{g,n})$ for
$i=0,1,\ldots, (3 g - 3+n)$ are defined as follows. Consider the
bundle $\omega_{g,n}(D) \, \to\, \C_{g,n}$ consisting of
$\omega_{g,n}$ twisted by the divisor $D = \sum_{i=1}^n D_i$, then
\[
\kappa_{(g,n),i} := \pi_*(c_1(\omega_{g,n}(D))^{i+1}).
\]
In particular, $\kappa_{(g,n),0} = 2g-2+n \in H^0(\M_{g,n})$ is the
negative of the Euler characteristic of a smooth curve of genus $g$
with $n$ points removed. We also have the equality $\omega_{g,n}(D) =
\LL_{(g,n+1),n+1}$ \cite{AC,HL}. Therefore
\begin{equation}
\label{psi:kappa}
\kappa_{(g,n),i} = \pi_*(\psi_{(g,n+1),n+1}^{i+1}).
\end{equation}

The tautological $\lambda$ classes are defined to be
\[
\lambda_{(g,n),i} := c_l (\pi_*\, \omega_{g,n}) 
\in H^{2i} (\M_{g,n}),
\]
where $l=1, \ldots, g$ because $\pi_*\, \omega_{g,n}$ is an orbifold
bundle of rank $g$. (There are no $\lambda$ classes in genus $0$ and
we define $\lambda_{(0,n),i} := 0$.) One can easily see that
$\lambda_{(g,n+1),i} = \pi^* \lambda_{(g,n),i}$. Therefore, all of the
$\lambda$ classes are pull backs of the $\lambda$ classes on
$\M_{1,1}$ and $\M_{g,0}$, $g\ge 2$. They can be expressed in terms of
the $\kappa$ classes, the $\psi$ classes, and the cohomology classes
lying at the boundary \cite{Mu}. In particular,
\begin{equation}
\label{lambda:kappa}
\kappa_{(g,n),1} = 12 \lambda_{(g,n),1} -\delta_{(g,n)} + 
\sum_{i=1}^n \psi_{(g,n),i},
\end{equation}
where $\delta_{(g,n)}$ is the fundamental class of $\M_{g,n} -
\MM_{g,n}$ \cite{Co,Fa2,Mu}. (This formula was brought to our attention by
E.~Getzler.)  

We will drop subscripts associated to the genus and the number of
punctures if there is no ambiguity.

\begin{nota}
Let $\SC_k$ be the set of infinite sequences of non-negative integers
$\mf = (m_k, m_{k+1}, m_{k+2}, \ldots)$ such that $m_i=0$ for all $i$
sufficiently large. We denote by $\delf_a$ the infinite sequence which
has only one non-zero entry 1 at the $a^{\text{th}}$ place. For $\mf =
(m_0, m_1, m_2, \ldots) \in \SC_0$ and $\tf = (t_0, t_1, t_2,
\ldots)$, a family of independent formal variables, we will use
notation of the type
\[
|\mf|:= \sum_{i\ge 0} i\, m_i, \quad |\!|\mf|\!|:= \sum_{i\ge 0} m_i,
\quad \mf! := \prod_{i\ge 0} m_i!, \quad \tf^\mf :=\prod_{i\ge 0}
t_i^{m_i}.
\]
We say that $\ll \le \mf$ if $l_i \le m_i$ for all $i$. If $\ll \le
\mf$ we let
\[
\binom{\mf}{\ll} := \prod_{i\ge 0} \binom{m_i}{l_i}
\]
We will use the same notation when $\mf \in \SC_1$. 
\end{nota}

\subsection{Generating Functions}

Witten \cite{W2} defined a generating function which incorporates all of
the information about the integrals of products of the $\psi$ classes.
In order to describe this function we need to introduce the following
notation. Let
\[
\la \tau_{d_1} \tau_{d_2} \ldots \tau_{d_n} \ra :=
\int_{\M_{g,n}} \psi_1^{d_1} \psi_2^{d_2} \ldots \psi_n^{d_n},
\]
where $g$ is determined by the equation $3g-3+n=d_1+d_2+\dots+d_n$. If
there exists no such $g$, then the left hand side is by definition
zero. In case we want to mention the genus explicitly we will write
$\la \tau_{d_1} \tau_{d_2} \ldots \tau_{d_n} \ra_g$. Note that this
expression is symmetric with respect to $d_1, d_2, \ldots, d_n$ since
$\psi_i$'s are interchanged under the action of the symmetric group
$S_n$. Therefore one can write it as $\la \tau_0^{m_0} \tau_1^{m_1}
\tau_2^{m_2} \ldots \ra$, where the set $\{ d_1,\ldots,d_n \}$
contains $m_0$ zeros, $m_1$ ones, etc. The generating function is
defined by
\[
F(t_0,t_1,t_2,\ldots) := \la \exp \sum_{j=0}^\infty t_j \tau_j \ra =
\sum_{ \{ m_i \} \in \SC_0} \prod_{i=0}^\infty \la \tau_0^{m_0}
\tau_1^{m_1} \tau_2^{m_2} \ldots \ra\, \frac{t_i^{m_i}}{m_i !}.
\]
We will also use the notation
\[
\sum_{\mf \in \SC_0} \la \tauf^\mf \ra \, \frac{\tf^\mf}{\mf!}
\]
for the last expression. Note that one can also write $F(\tf) =
\sum_{g=0}^\infty F_g(\tf)$, where $F_g(\tf) := \sum_{\mf \in \SC_0}
\la \tauf^\mf \ra_g\, \frac{\tf^\mf}{\mf!}$.

Witten conjectured in \cite{W} and Kontsevich proved in \cite{Ko} that
$F$ is the logarithm of a $\tau$-function in the KdV-hierarchy.

In \cite{KMZ} Kaufmann, Manin, and Zagier considered a similar
generating function for $\kappa$ classes. If one defines
\[
\la \kaf^\pf \ra_g =
\la \kappa_1^{p_1} \kappa_2^{p_2} \ldots \ra_g := 
\int_{\M_{g,n}} \kappa_1^{p_1} \kappa_2^{p_2} \ldots, 
\]
then their generating function is 
\[
K_g (x;\sf) = K_g(x;s_1, s_2, \ldots) :=
\sum_{\pf \in \SC_1} \la \kaf^\pf \ra_g\, \frac{x^{|\pf|}}{|\pf|!}\,
\frac{\sf^\pf}{\pf!}.
\]
Note that here it is important to indicate the genus. The number of
punctures $n$ is then determined from $3g-3+n= |\pf|$. 

We introduce the generating function $H$ which incorporates both of
the $\psi$ and $\kappa$ classes. We shall see that $F$ and $K$ enjoy
similar properties which arise because $H$ obeys those properties.

First we introduce the following notation. Let $\mf \in \SC_0$ and
$\pf \in \SC_1$. Define
\[
\la \tauf^\mf \kaf^\pf \ra :=
\int_{\M_{g,n}} \psi_1^{d_1} \psi_2^{d_2} \ldots \psi_n^{d_n}
            \kappa_1^{p_1} \kappa_2^{p_2} \ldots, 
\]
where the set $\{ d_1, \ldots, d_n \}$ contains $m_0$ zeros, $m_1$
ones, etc., and $(g,n)$ is determined by the equations $n=|\!| m
|\!|$, $3g-3+n= |\mf|+|\pf|$. If no such $g$ exists we set the
expression above to zero. As before, we write $\la
\tauf^\mf \kaf^\pf \ra_g$ when we want to fix $g$.

\begin{df}
Let $\tf= (t_0,t_1,\ldots)$ and $\sf=(s_1,s_2,\ldots)$ be independent
families of independent formal variables. We define
\[
H(\tf;\sf) := \sum_{\mf \in \SC_0,\pf \in \SC_1}
\la \tauf^\mf \kaf^\pf \ra\,
\frac{\tf^\mf}{\mf!}\, \frac{\sf^\pf}{\pf!}.
\]
\end{df}

One can split $H$ into the sum of $H_g$, $g=0,1,\ldots$. Each $H_g$
lies in a kernel of a certain scaling differential operator, i.e., it
satisfies the \emph{charge conservation equation}. The multiplication
of $H$ by $|\!| \mf |\!|$ is equivalent to applying the operator $\EE
:= \sum t_i \dd_i$, and by $|\mf|$ is equivalent to applying $\sum i\,
t_i \dd_i$. Therefore one has
\begin{equation}
\label{charge}
[ 3(1-g) + \sum_{i=0}^\infty (i-1)\, t_i \dd_i +
\sum_{i=1}^\infty i\, s_i d_i]\, H_g = 0,
\end{equation}
where $\dd_i:= \dd/\dd t_i$ and $d_i:= \dd/\dd s_i$. 

Clearly $H(\tf;\zef) = F(\tf)$. In order to relate
$H$ to $K$ one fixes a genus $g$, sets $t_1=t_2=\cdots=0$, and
$t_0=x$. The infinite sequence $\mf$ reduces to $m_0=n=|\pf|+3-3g$. It
follows that 
\[
H_g(x,\zef;\sf)= \sum_{\pf\in\SC_1} 
\la \kaf^\pf \ra_g\, \frac{x^{|\pf|+3-3g}}{(|\pf|+3-3g)!}\,
\frac{\sf^\pf}{\pf!}. 
\]
In this paper we are primarily interested in genus 0 and 1. Then
$K_0(x;\sf) = H_0'''(x,\zef;\sf)$ and $K_1(x;\sf) = H_1(x,\zef;\sf)$,
where the prime denotes the partial derivative with respect to $x$.

%	write kappas/psi in terms of graphs
% Last Modified 05/30/97 4pm by TK
\section{Presentation of the Tautological Classes via Graphs}
\label{graphs}

In this section we will mainly focus on $H_0$ and $H_1$, the
generating functions for the intersection numbers of the $\kappa$ and
$\psi$ classes in genus $0$ and genus $1$. First we show that $H_0$
satisfies a system of nonlinear differential equations. These
equations when, restricted to the $\psi$ classes, were first obtained
by Witten \cite{W2,W}. When restricted to the $\kappa$ classes in
genus zero, equivalent but much more complicated equations were
obtained by Kaufmann, Manin, and Zagier \cite{KMZ}, the difference
being accounted for by our simple, recursive presentation of powers of
$\kappa$ and $\psi$ classes in genus zero and one in terms of boundary
divisors. We prove that $H_0$ satisfies differential equations of the
same form as those obtained by Witten for just the $\psi$ classes. We
then obtain a system of differential equations relating $H_0$ and
$H_1$, and the explicit formula \eqref{our:rel}.

Here we present a geometric approach using the explicit presentation
of the $\psi$ and $\kappa$ classes in terms of the boundary strata.

\subsection{Basic Relation}
\label{basic}

In the beginning we want to state some general facts which hold for
all genera. Let $\pi: \C_{g,n} = \M_{g,n+1} \to \M_{g,n}$ be the
universal curve. To simplify the notation in this subsection, we denote 
$\psi_{(g,n),i}$ by $\psi_i$, $\psi_{(g,n+1),i}$ by $\psio_i$, and we
use the same convention for the $\kappa$ classes. Recall from
Sec.~\ref{kmz} that $D_i$ denotes the image of the $i^{\text{th}}$
canonical section $\sigma_i: \M_{g,n} \to \M_{g,n+1}$. We will denote
by the same letter its dual cohomology class in $H^2(\M_{g,n+1})$.

\begin{lm}
For each $k\ge 1$
\begin{equation}
\label{up}
\psio_i^a = \pi^* \psi_i^a + \sigma_{i*} \psi_i^{a-1}, \quad k\ge 1.
\end{equation}
\end{lm}

\begin{proof}
It is easy to show that the equation $\psio_i^a = \pi^* \psi_i^a +
\pi^* \psi_i^{a-1} D_i$ from \cite{W2} can be rewritten as
\begin{equation*}
\psio_i^a = \pi^* \psi_i^a + (-1)^{a-1} D_i^a, \quad k\ge 1. 
\end{equation*}
Applying the functor $\sigma_{i*} \sigma_i^*$, which is the
multiplication by $D_i$, to the above equation, and using that
$\sigma_i^* \psio_i=0$ one gets
\begin{equation*}
(-1)^{a} D_i^{a+1} = \sigma_{i*} \psi_i^{a}, \quad k\ge 0. 
\end{equation*}
These two equalities imply the lemma.
\end{proof}

\subsection{Explicit Presentation in Genus 0, 1}

In genus $0$ and genus $1$ equation \eqref{up} allows us to express
inductively all powers of $\psi$ classes, and therefore $\kappa$
classes, in terms of the boundary strata. Because the $\psi$ classes
are interchanged under the action of the symmetric group this is
enough to compute $\psi_{(g,n),1}$, $g=0,1$, and all its powers. In
the calculations below we will use the properties stated in
Sec.~\ref{kmz} regarding the pull backs and push forwards of the
cohomology classes represented by graphs.

Let us introduce the following notation for the rest of this section.
On graphs we denote by $\bullet$ vertices of genus $0$, and by $\circ$
vertices of genus $1$. We always assume that the $\psi$ classes are
associated to the marked point labeled $1$, and subsequently omit it
from the notation. We denote $\psi_{(1,n),1}$ by $\psi_{(n)}$, and we
denote $\psi_{(0,J),1}$ by $\phi_{(J)}$, where $J$ is a finite set,
$1\in J$. Similarly, we denote $\kappa_{(1,n),a}$ by $\kappa_{(n),a}$,
and $\kappa_{(0,J),a}$ by $\omega_{(J),a}$.

We also adopt the following convention. Let $\Gamma$ be a stable
graph. According to Sec.~\ref{kmz}, $\Gamma$ determines a canonical
finite quotient map from a product of moduli spaces to $\M_\Gamma$
provided certain choices have been made. We denote by $\rho_\Gamma$
the composition
\begin{equation}
\label{rho:gamma}
\rho_\Gamma: \prod_{v\in V(\Gamma)} \M_{g(v),n(v)} \longrightarrow
\M_\Gamma \longrightarrow \M_{g(\Gamma), S(\Gamma)},
\end{equation}
where the first arrow is the quotient morphism, and the second arrow
is the inclusion. 

Let $\gamma_v \in H^\bullet(\M_{g(v),n(v)})$. We denote
$\frac{1}{|\Aut(\Gamma)|} \rho_{\Gamma*} (\otimes_{v} \gamma_v)$ by
the picture of $\Gamma$ where each vertex $v$ is in addition labeled
by the cohomology class $\gamma_v$. We may omit the label of $v$ if
$\gamma(v)$ is the fundamental class of $\M_{g(v),n(v)}$. In
particular, the fundamental class of $\M_\Gamma$ (in the orbifold
sense) is represented by $\Gamma$ with all additional labels omitted.
The only possible ambiguity would arise if $\otimes_{v} \gamma_v$ were
not invariant under $\Aut (\Gamma)$, but this situation will not arise
in this paper.

In the pictures the dashed line with two arrows indicates the (sub)set
of the tails emanating from a particular vertex. If $\phi$ or $\omega$ labels
a vertex of a graph we will omit the subscript from the notation because
it is determined by the graph. (We also assume that $\phi$ is
associated to the marked point labeled by $1$.) The power $\phi^0$
represents the fundamental class.

\begin{prop}
\label{prop:zero:psi}
If $n\ge 4$, $a\ge 1$, then the following holds in $H^\bullet(\M_{0,n})$
\begin{equation}
\label{zero:psi}
\lift{1500}{$\phi_{(n)}^a = 
\displaystyle{\sum_{\substack{I\sqcup J = [n-2]\\ 1\in J}}}$}
\ \begin{picture}(0,0)%
\includegraphics{./pic/phi.pstex}%
\end{picture}%
\setlength{\unitlength}{0.00033300in}%
\begingroup\makeatletter\ifx\SetFigFont\undefined
% extract first six characters in \fmtname
\def\x#1#2#3#4#5#6#7\relax{\def\x{#1#2#3#4#5#6}}%
\expandafter\x\fmtname xxxxxx\relax \def\y{splain}%
\ifx\x\y   % LaTeX or SliTeX?
\gdef\SetFigFont#1#2#3{%
  \ifnum #1<17\tiny\else \ifnum #1<20\small\else
  \ifnum #1<24\normalsize\else \ifnum #1<29\large\else
  \ifnum #1<34\Large\else \ifnum #1<41\LARGE\else
     \huge\fi\fi\fi\fi\fi\fi
  \csname #3\endcsname}%
\else
\gdef\SetFigFont#1#2#3{\begingroup
  \count@#1\relax \ifnum 25<\count@\count@25\fi
  \def\x{\endgroup\@setsize\SetFigFont{#2pt}}%
  \expandafter\x
    \csname \romannumeral\the\count@ pt\expandafter\endcsname
    \csname @\romannumeral\the\count@ pt\endcsname
  \csname #3\endcsname}%
\fi
\fi\endgroup
\begin{picture}(3902,1479)(1276,-1873)
\put(2026,-586){\makebox(0,0)[lb]{\smash{\SetFigFont{8}{9.6}{rm}$n\!-\!1$}}}
\put(3001,-586){\makebox(0,0)[lb]{\smash{\SetFigFont{8}{9.6}{rm}$n$}}}
\put(5178,-1263){\makebox(0,0)[lb]{\smash{\SetFigFont{8}{9.6}{rm}$J$}}}
\put(1276,-1261){\makebox(0,0)[lb]{\smash{\SetFigFont{8}{9.6}{rm}$I$}}}
\put(4951,-586){\makebox(0,0)[lb]{\smash{\SetFigFont{8}{9.6}{rm}$1$}}}
\put(3301,-1711){\makebox(0,0)[lb]{\smash{\SetFigFont{8}{9.6}{rm}$\phi^{a-1}$}}}
\end{picture}

%\phi_{(n)}^a = \sum_{\substack{I\sqcup J = [n-2]\\ 1\in J}}
%A_I \otimes \phi_{(J)}^{a-1}.
\end{equation}
\end{prop}

\begin{rem}
The class $\phi_{(n)}^a$ is invariant under the subgroup $S_{n-1}
\subset S_n$ whose elements fix $1$. Therefore instead of $n-1,n$ we
can choose any two labels $a,b$, $2\le a<b \le n$ to be distinguished.
\end{rem}

\begin{proof}
The statement is true when $n=4$ and $a=1$ because $\int_{\M_{0,4}}
\psi_1 = 1$. We shall first prove by induction that the statement is
true for all $n\ge 4$ and $a=1$. Let us consider the projection $\pi:
\M_{0,n+1} \to \M_{0,n}$ which ``forgets'' the $n+1^{\text{st}}$
marked point.
%If $I$ is a finite set we denote by $I^+$ the set $I \sqcup \{ n+1 \}$.
By the induction hypothesis we assume that the statement is true for
some $n$ and $a=1$. Applying $\pi^*$ one gets:
\begin{align*}
\lift{1500}{$\phi_{(n+1)}\ - $} \quad
\ \begin{picture}(0,0)%
\includegraphics{./pic/d1.pstex}%
\end{picture}%
\setlength{\unitlength}{0.00033300in}%
\begingroup\makeatletter\ifx\SetFigFont\undefined
% extract first six characters in \fmtname
\def\x#1#2#3#4#5#6#7\relax{\def\x{#1#2#3#4#5#6}}%
\expandafter\x\fmtname xxxxxx\relax \def\y{splain}%
\ifx\x\y   % LaTeX or SliTeX?
\gdef\SetFigFont#1#2#3{%
  \ifnum #1<17\tiny\else \ifnum #1<20\small\else
  \ifnum #1<24\normalsize\else \ifnum #1<29\large\else
  \ifnum #1<34\Large\else \ifnum #1<41\LARGE\else
     \huge\fi\fi\fi\fi\fi\fi
  \csname #3\endcsname}%
\else
\gdef\SetFigFont#1#2#3{\begingroup
  \count@#1\relax \ifnum 25<\count@\count@25\fi
  \def\x{\endgroup\@setsize\SetFigFont{#2pt}}%
  \expandafter\x
    \csname \romannumeral\the\count@ pt\expandafter\endcsname
    \csname @\romannumeral\the\count@ pt\endcsname
  \csname #3\endcsname}%
\fi
\fi\endgroup
\begin{picture}(4425,1515)(526,-1909)
\put(4951,-586){\makebox(0,0)[lb]{\smash{\SetFigFont{8}{9.6}{rm}$1$}}}
\put(4951,-1861){\makebox(0,0)[lb]{\smash{\SetFigFont{8}{9.6}{rm}$n\!+\!1$}}}
\put(526,-1336){\makebox(0,0)[lb]{\smash{\SetFigFont{8}{9.6}{rm}$[n\!-\!2]$}}}
\put(3001,-586){\makebox(0,0)[lb]{\smash{\SetFigFont{8}{9.6}{rm}$n$}}}
\put(2026,-586){\makebox(0,0)[lb]{\smash{\SetFigFont{8}{9.6}{rm}$n\!-\!1$}}}
\end{picture}
 &\\
\lift{1500}
{$= \displaystyle{\sum_{\substack{I\sqcup J = [n-2]\\ 1\in J}}}
\bigg( $}
\ \begin{picture}(0,0)%
\includegraphics{./pic/proof1.pstex}%
\end{picture}%
\setlength{\unitlength}{0.00033300in}%
\begingroup\makeatletter\ifx\SetFigFont\undefined
% extract first six characters in \fmtname
\def\x#1#2#3#4#5#6#7\relax{\def\x{#1#2#3#4#5#6}}%
\expandafter\x\fmtname xxxxxx\relax \def\y{splain}%
\ifx\x\y   % LaTeX or SliTeX?
\gdef\SetFigFont#1#2#3{%
  \ifnum #1<17\tiny\else \ifnum #1<20\small\else
  \ifnum #1<24\normalsize\else \ifnum #1<29\large\else
  \ifnum #1<34\Large\else \ifnum #1<41\LARGE\else
     \huge\fi\fi\fi\fi\fi\fi
  \csname #3\endcsname}%
\else
\gdef\SetFigFont#1#2#3{\begingroup
  \count@#1\relax \ifnum 25<\count@\count@25\fi
  \def\x{\endgroup\@setsize\SetFigFont{#2pt}}%
  \expandafter\x
    \csname \romannumeral\the\count@ pt\expandafter\endcsname
    \csname @\romannumeral\the\count@ pt\endcsname
  \csname #3\endcsname}%
\fi
\fi\endgroup
\begin{picture}(3902,1479)(1276,-1873)
\put(5178,-1263){\makebox(0,0)[lb]{\smash{\SetFigFont{8}{9.6}{rm}$J$}}}
\put(1276,-1261){\makebox(0,0)[lb]{\smash{\SetFigFont{8}{9.6}{rm}$I$}}}
\put(4951,-586){\makebox(0,0)[lb]{\smash{\SetFigFont{8}{9.6}{rm}$1$}}}
\put(3376,-586){\makebox(0,0)[lb]{\smash{\SetFigFont{8}{9.6}{rm}$n\!+\!1$}}}
\put(2926,-586){\makebox(0,0)[lb]{\smash{\SetFigFont{8}{9.6}{rm}$n$}}}
\put(2026,-586){\makebox(0,0)[lb]{\smash{\SetFigFont{8}{9.6}{rm}$n\!-\!1$}}}
\end{picture}
 &
\lift{1500}{$+$} 
\quad \begin{picture}(0,0)%
\includegraphics{./pic/proof2.pstex}%
\end{picture}%
\setlength{\unitlength}{0.00033300in}%
\begingroup\makeatletter\ifx\SetFigFont\undefined
% extract first six characters in \fmtname
\def\x#1#2#3#4#5#6#7\relax{\def\x{#1#2#3#4#5#6}}%
\expandafter\x\fmtname xxxxxx\relax \def\y{splain}%
\ifx\x\y   % LaTeX or SliTeX?
\gdef\SetFigFont#1#2#3{%
  \ifnum #1<17\tiny\else \ifnum #1<20\small\else
  \ifnum #1<24\normalsize\else \ifnum #1<29\large\else
  \ifnum #1<34\Large\else \ifnum #1<41\LARGE\else
     \huge\fi\fi\fi\fi\fi\fi
  \csname #3\endcsname}%
\else
\gdef\SetFigFont#1#2#3{\begingroup
  \count@#1\relax \ifnum 25<\count@\count@25\fi
  \def\x{\endgroup\@setsize\SetFigFont{#2pt}}%
  \expandafter\x
    \csname \romannumeral\the\count@ pt\expandafter\endcsname
    \csname @\romannumeral\the\count@ pt\endcsname
  \csname #3\endcsname}%
\fi
\fi\endgroup
\begin{picture}(3902,1479)(1276,-1873)
\put(5178,-1263){\makebox(0,0)[lb]{\smash{\SetFigFont{8}{9.6}{rm}$J$}}}
\put(1276,-1261){\makebox(0,0)[lb]{\smash{\SetFigFont{8}{9.6}{rm}$I$}}}
\put(4951,-586){\makebox(0,0)[lb]{\smash{\SetFigFont{8}{9.6}{rm}$1$}}}
\put(3676,-586){\makebox(0,0)[lb]{\smash{\SetFigFont{8}{9.6}{rm}$n\!+\!1$}}}
\put(3001,-586){\makebox(0,0)[lb]{\smash{\SetFigFont{8}{9.6}{rm}$n$}}}
\put(2026,-586){\makebox(0,0)[lb]{\smash{\SetFigFont{8}{9.6}{rm}$n\!-\!1$}}}
\end{picture}
 \lift{1500}{$\!\!\!\! \bigg) $}
%\phi_{(n+1)} - D_1 = \sum_{\substack{I\sqcup J = [n-2]\\ 1\in J}}
%(A_{I^+} \otimes \phi_{(J)}^0 + A_I \otimes \phi_{(J^+)}^0).
\end{align*}
The left hand side is equal to $\phi_{(n+1)} - D_1 = \pi^*
\phi_{(n)}$. Moving $D_1$ to the right hand side and relabeling the
tails marked $n-1, n$ by $n,n+1$ respectively one gets the statement
in case of $n+1$, $a=1$.

%Therefore
%\[
%\phi_{(n+1)} = \sum_{\substack{I\sqcup J = [n-1]\\ 1\in J}}
%A_I \otimes \phi_{(J)}^0. 
%\]
%Here $D_1$ corresponds to $J=\{ 1,n-1 \} = \{ 1 \}^+$. 

Now assume that $a\ge 2$. If $n=4$, then the statement of the
proposition is trivially satisfied because all terms vanish by
dimensional considerations. We assume that the statement is true for a
pair $(n,a)$ and all pairs $(n',a')$, where $a'<a$, and prove it for
$(n+1,a)$. Applying $\pi^*$ one gets
\begin{align*}
\lift{1500}{$\phi_{(n+1)}^a - \sigma_{1*} \phi_{(n)}^{a-1} = 
\displaystyle{\sum_{\substack{I\sqcup J = [n-2]\\ 1\in J}}} 
\bigg( $}
\ \begin{picture}(0,0)%
\includegraphics{./pic/proof3.pstex}%
\end{picture}%
\setlength{\unitlength}{0.00033300in}%
\begingroup\makeatletter\ifx\SetFigFont\undefined
% extract first six characters in \fmtname
\def\x#1#2#3#4#5#6#7\relax{\def\x{#1#2#3#4#5#6}}%
\expandafter\x\fmtname xxxxxx\relax \def\y{splain}%
\ifx\x\y   % LaTeX or SliTeX?
\gdef\SetFigFont#1#2#3{%
  \ifnum #1<17\tiny\else \ifnum #1<20\small\else
  \ifnum #1<24\normalsize\else \ifnum #1<29\large\else
  \ifnum #1<34\Large\else \ifnum #1<41\LARGE\else
     \huge\fi\fi\fi\fi\fi\fi
  \csname #3\endcsname}%
\else
\gdef\SetFigFont#1#2#3{\begingroup
  \count@#1\relax \ifnum 25<\count@\count@25\fi
  \def\x{\endgroup\@setsize\SetFigFont{#2pt}}%
  \expandafter\x
    \csname \romannumeral\the\count@ pt\expandafter\endcsname
    \csname @\romannumeral\the\count@ pt\endcsname
  \csname #3\endcsname}%
\fi
\fi\endgroup
\begin{picture}(3902,1479)(1276,-1873)
\put(5178,-1263){\makebox(0,0)[lb]{\smash{\SetFigFont{8}{9.6}{rm}$J$}}}
\put(1276,-1261){\makebox(0,0)[lb]{\smash{\SetFigFont{8}{9.6}{rm}$I$}}}
\put(4951,-586){\makebox(0,0)[lb]{\smash{\SetFigFont{8}{9.6}{rm}$1$}}}
\put(3376,-586){\makebox(0,0)[lb]{\smash{\SetFigFont{8}{9.6}{rm}$n\!+\!1$}}}
\put(2926,-586){\makebox(0,0)[lb]{\smash{\SetFigFont{8}{9.6}{rm}$n$}}}
\put(2026,-586){\makebox(0,0)[lb]{\smash{\SetFigFont{8}{9.6}{rm}$n\!-\!1$}}}
\put(3301,-1711){\makebox(0,0)[lb]{\smash{\SetFigFont{8}{9.6}{rm}$\phi^{a-1}$}}}
\end{picture}
 & \\
\lift{1500}{$+$} \quad \begin{picture}(0,0)%
\includegraphics{./pic/proof4.pstex}%
\end{picture}%
\setlength{\unitlength}{0.00033300in}%
\begingroup\makeatletter\ifx\SetFigFont\undefined
% extract first six characters in \fmtname
\def\x#1#2#3#4#5#6#7\relax{\def\x{#1#2#3#4#5#6}}%
\expandafter\x\fmtname xxxxxx\relax \def\y{splain}%
\ifx\x\y   % LaTeX or SliTeX?
\gdef\SetFigFont#1#2#3{%
  \ifnum #1<17\tiny\else \ifnum #1<20\small\else
  \ifnum #1<24\normalsize\else \ifnum #1<29\large\else
  \ifnum #1<34\Large\else \ifnum #1<41\LARGE\else
     \huge\fi\fi\fi\fi\fi\fi
  \csname #3\endcsname}%
\else
\gdef\SetFigFont#1#2#3{\begingroup
  \count@#1\relax \ifnum 25<\count@\count@25\fi
  \def\x{\endgroup\@setsize\SetFigFont{#2pt}}%
  \expandafter\x
    \csname \romannumeral\the\count@ pt\expandafter\endcsname
    \csname @\romannumeral\the\count@ pt\endcsname
  \csname #3\endcsname}%
\fi
\fi\endgroup
\begin{picture}(3902,1479)(1276,-1873)
\put(5178,-1263){\makebox(0,0)[lb]{\smash{\SetFigFont{8}{9.6}{rm}$J$}}}
\put(1276,-1261){\makebox(0,0)[lb]{\smash{\SetFigFont{8}{9.6}{rm}$I$}}}
\put(4951,-586){\makebox(0,0)[lb]{\smash{\SetFigFont{8}{9.6}{rm}$1$}}}
\put(3676,-586){\makebox(0,0)[lb]{\smash{\SetFigFont{8}{9.6}{rm}$n\!+\!1$}}}
\put(3001,-586){\makebox(0,0)[lb]{\smash{\SetFigFont{8}{9.6}{rm}$n$}}}
\put(2026,-586){\makebox(0,0)[lb]{\smash{\SetFigFont{8}{9.6}{rm}$n\!-\!1$}}}
\put(3301,-1711){\makebox(0,0)[lb]{\smash{\SetFigFont{8}{9.6}{rm}$\phi^{a-1}$}}}
\end{picture}
 
\lift{1500}{$-\ \sigma_{1*} \big( \ $} \begin{picture}(0,0)%
\includegraphics{./pic/sigma.pstex}%
\end{picture}%
\setlength{\unitlength}{0.00033300in}%
\begingroup\makeatletter\ifx\SetFigFont\undefined
% extract first six characters in \fmtname
\def\x#1#2#3#4#5#6#7\relax{\def\x{#1#2#3#4#5#6}}%
\expandafter\x\fmtname xxxxxx\relax \def\y{splain}%
\ifx\x\y   % LaTeX or SliTeX?
\gdef\SetFigFont#1#2#3{%
  \ifnum #1<17\tiny\else \ifnum #1<20\small\else
  \ifnum #1<24\normalsize\else \ifnum #1<29\large\else
  \ifnum #1<34\Large\else \ifnum #1<41\LARGE\else
     \huge\fi\fi\fi\fi\fi\fi
  \csname #3\endcsname}%
\else
\gdef\SetFigFont#1#2#3{\begingroup
  \count@#1\relax \ifnum 25<\count@\count@25\fi
  \def\x{\endgroup\@setsize\SetFigFont{#2pt}}%
  \expandafter\x
    \csname \romannumeral\the\count@ pt\expandafter\endcsname
    \csname @\romannumeral\the\count@ pt\endcsname
  \csname #3\endcsname}%
\fi
\fi\endgroup
\begin{picture}(3902,1479)(1276,-1873)
\put(2026,-586){\makebox(0,0)[lb]{\smash{\SetFigFont{8}{9.6}{rm}$n\!-\!1$}}}
\put(3001,-586){\makebox(0,0)[lb]{\smash{\SetFigFont{8}{9.6}{rm}$n$}}}
\put(5178,-1263){\makebox(0,0)[lb]{\smash{\SetFigFont{8}{9.6}{rm}$J$}}}
\put(1276,-1261){\makebox(0,0)[lb]{\smash{\SetFigFont{8}{9.6}{rm}$I$}}}
\put(4951,-586){\makebox(0,0)[lb]{\smash{\SetFigFont{8}{9.6}{rm}$1$}}}
\put(3301,-1711){\makebox(0,0)[lb]{\smash{\SetFigFont{8}{9.6}{rm}$\phi^{a-2}$}}}
\end{picture}
 &
\lift{1500}{$\!\!\!\! \big) \bigg) $}
%\phi_{(n+1)}^a - \sigma_{1*} \phi_{(n)}^{a-1} = 
%\sum_{\substack{I\sqcup J = [n-2], 1\in J}}
%(A_{I^+} \otimes \phi_{(J)}^{a-1} + A_I \otimes 
%\phi_{(J^+)}^{a-1} - A_I \otimes \sigma_{1*} \phi_{(J)}^{a-2} ).
\end{align*}

According to the induction hypothesis the terms with $\sigma_{1*}$ on
the left hand side and the right hand side cancel each other. Thus, we
get the statement of the proposition for the pair $(n+1,a)$. 
\end{proof}

Let $\pi_1$ be the morphism $\M_{0,n+1} \to \M_{0,n}$ forgetting the
first marked point. Applying it to \eqref{zero:psi} with $a+1$, using
\eqref{psi:kappa}, and renumbering the labels $\{ 2,\ldots,n+1 \}$ by
the elements of $[n]$ we get the following

\begin{crl}
If $n\ge 4$, $a\ge1$, then the following holds in $H^\bullet(\M_{0,n})$
\begin{equation}
\label{zero:kappa}
\lift{1500}{$\omega_{(n),a} =
\displaystyle{\sum_{I\sqcup J = [n-2]}}$}
\ \begin{picture}(0,0)%
\includegraphics{./pic/omega.pstex}%
\end{picture}%
\setlength{\unitlength}{0.00033300in}%
\begingroup\makeatletter\ifx\SetFigFont\undefined
% extract first six characters in \fmtname
\def\x#1#2#3#4#5#6#7\relax{\def\x{#1#2#3#4#5#6}}%
\expandafter\x\fmtname xxxxxx\relax \def\y{splain}%
\ifx\x\y   % LaTeX or SliTeX?
\gdef\SetFigFont#1#2#3{%
  \ifnum #1<17\tiny\else \ifnum #1<20\small\else
  \ifnum #1<24\normalsize\else \ifnum #1<29\large\else
  \ifnum #1<34\Large\else \ifnum #1<41\LARGE\else
     \huge\fi\fi\fi\fi\fi\fi
  \csname #3\endcsname}%
\else
\gdef\SetFigFont#1#2#3{\begingroup
  \count@#1\relax \ifnum 25<\count@\count@25\fi
  \def\x{\endgroup\@setsize\SetFigFont{#2pt}}%
  \expandafter\x
    \csname \romannumeral\the\count@ pt\expandafter\endcsname
    \csname @\romannumeral\the\count@ pt\endcsname
  \csname #3\endcsname}%
\fi
\fi\endgroup
\begin{picture}(3902,1515)(1276,-1873)
\put(2026,-586){\makebox(0,0)[lb]{\smash{\SetFigFont{8}{9.6}{rm}$n\!-\!1$}}}
\put(3001,-586){\makebox(0,0)[lb]{\smash{\SetFigFont{8}{9.6}{rm}$n$}}}
\put(5178,-1263){\makebox(0,0)[lb]{\smash{\SetFigFont{8}{9.6}{rm}$J$}}}
\put(1276,-1261){\makebox(0,0)[lb]{\smash{\SetFigFont{8}{9.6}{rm}$I$}}}
\put(4951,-586){\makebox(0,0)[lb]{\smash{\SetFigFont{8}{9.6}{rm}$\phantom{1}$}}}
\put(3451,-1711){\makebox(0,0)[lb]{\smash{\SetFigFont{8}{9.6}{rm}$\omega_{a-1}$}}}
\end{picture}
 \qed
%\omega_{(n),a} = \sum_{I\sqcup J = [n-2]} 
%A_I \otimes \omega_{(J),a-1}. \qed
\end{equation}
\end{crl}

In \eqref{zero:kappa} if $a=1$ one should use that $\omega_{(J\sqcup
*),0}$ associated to the right vertex is equal to $|J|-1$ times the
fundamental class (cf.~Sec.~\ref{kmz}). This agrees with $\pi_{1*}$ of
\eqref{zero:psi} when $a=1$.

Now we establish a relation between the $\psi$ and $\kappa$ classes in
genus $0$ and $1$. The proof is virtually identical to that of
Prop.~\ref{prop:zero:psi}, and we will not reproduce it. 
%Suppose that $I\sqcup J=[n]$. Then the following graphs
%\[
%figure \Gamma_3 \Gamma_4, \Gamma_5, \Gamma_6
%\]
%determine the elements $\phi_{(n)}^{'a}, 
%\omega'_{(n),a}$ and $B_I \otimes \phi_{(J)}^a, B_I \otimes
%\omega_{(J),a}$ in $H^\bullet(\M_{1,n})$.
%This is shown in \cite[VI.4]{DR} that $\psi_{(1)} = \frac{1}{12}
%\phi_{(1)}^{'0} $. Using this and \eqref{up} one gets

It is shown in \cite[VI.4]{DR} that 
\begin{equation*}
\lift{850}{$\psi_{(1)} = \displaystyle{\frac{1}{12}}$}
\ \begin{picture}(0,0)%
\includegraphics{./pic/psi2.pstex}%
\end{picture}%
\setlength{\unitlength}{0.00033300in}%
\begingroup\makeatletter\ifx\SetFigFont\undefined
% extract first six characters in \fmtname
\def\x#1#2#3#4#5#6#7\relax{\def\x{#1#2#3#4#5#6}}%
\expandafter\x\fmtname xxxxxx\relax \def\y{splain}%
\ifx\x\y   % LaTeX or SliTeX?
\gdef\SetFigFont#1#2#3{%
  \ifnum #1<17\tiny\else \ifnum #1<20\small\else
  \ifnum #1<24\normalsize\else \ifnum #1<29\large\else
  \ifnum #1<34\Large\else \ifnum #1<41\LARGE\else
     \huge\fi\fi\fi\fi\fi\fi
  \csname #3\endcsname}%
\else
\gdef\SetFigFont#1#2#3{\begingroup
  \count@#1\relax \ifnum 25<\count@\count@25\fi
  \def\x{\endgroup\@setsize\SetFigFont{#2pt}}%
  \expandafter\x
    \csname \romannumeral\the\count@ pt\expandafter\endcsname
    \csname @\romannumeral\the\count@ pt\endcsname
  \csname #3\endcsname}%
\fi
\fi\endgroup
\begin{picture}(2484,853)(1492,-1687)
\put(3976,-1336){\makebox(0,0)[lb]{\smash{\SetFigFont{8}{9.6}{rm}$1$}}}
\end{picture}

\end{equation*}
Note that we take the coefficient is $\frac{1}{12}$ rather than
$\frac{1}{24}$ due to the non-trivial automorphism of the graph. Using
this and \eqref{up} one obtains the following
\begin{prop}
\label{prop:one:psi}
If $n\ge 1$, $a\ge 1$, then the following holds in $H^\bullet(\M_{1,n})$
\begin{equation}
\label{one:psi}
\lift{1500}{$\psi_{(n)}^a = \displaystyle{\frac{1}{12}}$}
\ \begin{picture}(0,0)%
\includegraphics{./pic/loopsi.pstex}%
\end{picture}%
\setlength{\unitlength}{0.00033300in}%
\begingroup\makeatletter\ifx\SetFigFont\undefined
% extract first six characters in \fmtname
\def\x#1#2#3#4#5#6#7\relax{\def\x{#1#2#3#4#5#6}}%
\expandafter\x\fmtname xxxxxx\relax \def\y{splain}%
\ifx\x\y   % LaTeX or SliTeX?
\gdef\SetFigFont#1#2#3{%
  \ifnum #1<17\tiny\else \ifnum #1<20\small\else
  \ifnum #1<24\normalsize\else \ifnum #1<29\large\else
  \ifnum #1<34\Large\else \ifnum #1<41\LARGE\else
     \huge\fi\fi\fi\fi\fi\fi
  \csname #3\endcsname}%
\else
\gdef\SetFigFont#1#2#3{\begingroup
  \count@#1\relax \ifnum 25<\count@\count@25\fi
  \def\x{\endgroup\@setsize\SetFigFont{#2pt}}%
  \expandafter\x
    \csname \romannumeral\the\count@ pt\expandafter\endcsname
    \csname @\romannumeral\the\count@ pt\endcsname
  \csname #3\endcsname}%
\fi
\fi\endgroup
\begin{picture}(2487,1480)(1492,-1874)
\put(3979,-1264){\makebox(0,0)[lb]{\smash{\SetFigFont{8}{9.6}{rm}$[n]$}}}
\put(3751,-586){\makebox(0,0)[lb]{\smash{\SetFigFont{8}{9.6}{rm}$1$}}}
\put(2251,-1786){\makebox(0,0)[lb]{\smash{\SetFigFont{8}{9.6}{rm}$\phi^{a-1}$}}}
\end{picture}
  
\lift{1500}{$\displaystyle{+ \sum_{\substack{I\sqcup J=[n]\\ 1\in J} }}$}
\ \begin{picture}(0,0)%
\includegraphics{./pic/psi.pstex}%
\end{picture}%
\setlength{\unitlength}{0.00033300in}%
\begingroup\makeatletter\ifx\SetFigFont\undefined
% extract first six characters in \fmtname
\def\x#1#2#3#4#5#6#7\relax{\def\x{#1#2#3#4#5#6}}%
\expandafter\x\fmtname xxxxxx\relax \def\y{splain}%
\ifx\x\y   % LaTeX or SliTeX?
\gdef\SetFigFont#1#2#3{%
  \ifnum #1<17\tiny\else \ifnum #1<20\small\else
  \ifnum #1<24\normalsize\else \ifnum #1<29\large\else
  \ifnum #1<34\Large\else \ifnum #1<41\LARGE\else
     \huge\fi\fi\fi\fi\fi\fi
  \csname #3\endcsname}%
\else
\gdef\SetFigFont#1#2#3{\begingroup
  \count@#1\relax \ifnum 25<\count@\count@25\fi
  \def\x{\endgroup\@setsize\SetFigFont{#2pt}}%
  \expandafter\x
    \csname \romannumeral\the\count@ pt\expandafter\endcsname
    \csname @\romannumeral\the\count@ pt\endcsname
  \csname #3\endcsname}%
\fi
\fi\endgroup
\begin{picture}(3902,1479)(1276,-1873)
\put(5178,-1263){\makebox(0,0)[lb]{\smash{\SetFigFont{8}{9.6}{rm}$J$}}}
\put(1276,-1261){\makebox(0,0)[lb]{\smash{\SetFigFont{8}{9.6}{rm}$I$}}}
\put(4951,-586){\makebox(0,0)[lb]{\smash{\SetFigFont{8}{9.6}{rm}$1$}}}
\put(3301,-1711){\makebox(0,0)[lb]{\smash{\SetFigFont{8}{9.6}{rm}$\phi^{a-1}$}}}
\end{picture}
 \qed
%\psi_{(n)}^a = \frac{1}{12}  \phi_{(n)}^{'a-1} +
%\sum_{\substack{I\sqcup J=[n]\\ 1\in J} }
%B_I \otimes \phi_{(J)}^{a-1}. \qed
\end{equation}
\end{prop}

Pushing down the above along $\pi_1: \M_{1,n+1} \to \M_{1,n}$, and
renumbering the labels one gets

\begin{crl}
If $n\ge 1$, $a\ge 1$, then the following holds in $H^\bullet(\M_{1,n})$
\begin{equation}
\label{one:kappa}
\lift{1500}{$\kappa_{(n),a}= \displaystyle{\frac{1}{12}}$}
\ \begin{picture}(0,0)%
\includegraphics{./pic/loopka.pstex}%
\end{picture}%
\setlength{\unitlength}{0.00033300in}%
\begingroup\makeatletter\ifx\SetFigFont\undefined
% extract first six characters in \fmtname
\def\x#1#2#3#4#5#6#7\relax{\def\x{#1#2#3#4#5#6}}%
\expandafter\x\fmtname xxxxxx\relax \def\y{splain}%
\ifx\x\y   % LaTeX or SliTeX?
\gdef\SetFigFont#1#2#3{%
  \ifnum #1<17\tiny\else \ifnum #1<20\small\else
  \ifnum #1<24\normalsize\else \ifnum #1<29\large\else
  \ifnum #1<34\Large\else \ifnum #1<41\LARGE\else
     \huge\fi\fi\fi\fi\fi\fi
  \csname #3\endcsname}%
\else
\gdef\SetFigFont#1#2#3{\begingroup
  \count@#1\relax \ifnum 25<\count@\count@25\fi
  \def\x{\endgroup\@setsize\SetFigFont{#2pt}}%
  \expandafter\x
    \csname \romannumeral\the\count@ pt\expandafter\endcsname
    \csname @\romannumeral\the\count@ pt\endcsname
  \csname #3\endcsname}%
\fi
\fi\endgroup
\begin{picture}(2487,1516)(1492,-1874)
\put(3979,-1264){\makebox(0,0)[lb]{\smash{\SetFigFont{8}{9.6}{rm}$[n]$}}}
\put(3751,-586){\makebox(0,0)[lb]{\smash{\SetFigFont{8}{9.6}{rm}$\phantom{1}$}}}
\put(2401,-1786){\makebox(0,0)[lb]{\smash{\SetFigFont{8}{9.6}{rm}$\omega_{a-1}$}}}
\end{picture}
  
\lift{1500}{$\displaystyle{+ \sum_{I\sqcup J=[n]}}$}
\ \begin{picture}(0,0)%
\includegraphics{./pic/kappa.pstex}%
\end{picture}%
\setlength{\unitlength}{0.00033300in}%
\begingroup\makeatletter\ifx\SetFigFont\undefined
% extract first six characters in \fmtname
\def\x#1#2#3#4#5#6#7\relax{\def\x{#1#2#3#4#5#6}}%
\expandafter\x\fmtname xxxxxx\relax \def\y{splain}%
\ifx\x\y   % LaTeX or SliTeX?
\gdef\SetFigFont#1#2#3{%
  \ifnum #1<17\tiny\else \ifnum #1<20\small\else
  \ifnum #1<24\normalsize\else \ifnum #1<29\large\else
  \ifnum #1<34\Large\else \ifnum #1<41\LARGE\else
     \huge\fi\fi\fi\fi\fi\fi
  \csname #3\endcsname}%
\else
\gdef\SetFigFont#1#2#3{\begingroup
  \count@#1\relax \ifnum 25<\count@\count@25\fi
  \def\x{\endgroup\@setsize\SetFigFont{#2pt}}%
  \expandafter\x
    \csname \romannumeral\the\count@ pt\expandafter\endcsname
    \csname @\romannumeral\the\count@ pt\endcsname
  \csname #3\endcsname}%
\fi
\fi\endgroup
\begin{picture}(3902,1515)(1276,-1873)
\put(5178,-1263){\makebox(0,0)[lb]{\smash{\SetFigFont{8}{9.6}{rm}$J$}}}
\put(1276,-1261){\makebox(0,0)[lb]{\smash{\SetFigFont{8}{9.6}{rm}$I$}}}
\put(3451,-1711){\makebox(0,0)[lb]{\smash{\SetFigFont{8}{9.6}{rm}$\omega_{a-1}$}}}
\put(4951,-586){\makebox(0,0)[lb]{\smash{\SetFigFont{8}{9.6}{rm}$\phantom{1}$}}}
\end{picture}
 \qed
%\kappa_{(n),a} = \frac{1}{12} \omega'_{(n),a-1} +
%\sum_{I\sqcup J=[n]} B_I \otimes \omega_{(J),a-1}. \qed
\end{equation}
\end{crl}

\subsection{Recursion Relations and Differential Equations}

Now we derive the corresponding recursion relations and differential
equations for the intersection numbers of the products of the $\psi$
and $\kappa$ classes using the explicit graph presentations above. In
order to obtain the recursion relations we use a method from
\cite{KMZ} to integrate the product of the $\psi$ and $\kappa$ classes
over the Poincar\'e dual of a chosen $\psi$ or $\kappa$ class. In
order to do this we need to know how the tautological classes restrict
to the strata of the natural stratification. The restriction of a
$\psi$ class to a boundary stratum is obvious. 

In order to restrict products of the $\kappa$ classes we use Lemma 1.3
from \cite{KMZ} where the authors show the following restriction
property for the $\kappa$ classes. (They show it in the case of genus
$0$, but their proof is in fact valid for all genera.) Let
$(\Gamma,g,\mu)$ be a stable graph, $\rho_\Gamma$ is the corresponding
morphism defined by \eqref{rho:gamma}, and $\kaf^\pf$ is a product of
the $\kappa$ classes on $\M_{g(\Gamma),S(\Gamma)}$. Then
\[
\int_{\prod_{v\in V(\Gamma)} \M_{g(v),n(v)}}
\frac{\rho_\Gamma^* (\kaf^\pf)}{\pf !} = 
\sum_{\substack{\pf^v:\,v\in V(\Gamma)\\ \sum \pf^v = \pf}}\;
\prod_{v\in V(\Gamma)} 
\frac{\la \kaf^{\pf^v} \ra_{g(v)}} {\pf^v !}
\]
The argument uses a fact proved in \cite{AC} that the collection
$\kappa_{(g,n),a}$ for each fixed $a$ forms a \emph{logarithmic
cohomological field theory} (cf. Sec.~\ref{cft}), i.e., the $\kappa$
classes satisfy the relation

\begin{equation}
\label{eq:lcft}
\rho_\Gamma^* (\kappa_a) = \sum_{v\in V(\Gamma)} \kappa_{g(v),n(v)}.
\end{equation}

We start with genus $0$. Recall that $H(\tf;\sf)$ is the generating
function incorporating the intersection numbers for the products of
the $\psi$ and $\kappa$ classes defined in Sec.~\ref{kmz}, $\dd_a$,
$d_a$ are partial derivative with respect to $t_a$, $a\ge 0$, $s_a$,
$a\ge 1$, and $\EE = \sum_{i=0}^\infty t_i\dd_i$.

\begin{thm}
\label{H0}
For each $\mf \in \SC_1$, $\pf \in \SC$, $k,l \ge 0$, and $a\ge 1$ one
has
\begin{align*}
\la \tauf^{\mf + \delf_k + \delf_l + \delf_a} \kaf^\pf \ra_0 & \\
= \sum_{\substack{\mf'+\mf''=\mf\\ \pf' + \pf''= \pf} }
& \binom{\mf}{\mf'} \binom{\pf}{\pf'} 
\la \tauf^{\mf' + \delf_k + \delf_l + \delf_0} \kaf^{\pf'} \ra_0
\la \tauf^{\mf'' + \delf_{a-1} + \delf_0} \kaf^{\pf''} \ra_0, \\
\la \tauf^{\mf + \delf_k + \delf_l} \kaf^{\pf + \delf_a} \ra_0 & \\
= \sum_{\substack{\mf'+\mf''=\mf\\ \pf' + \pf''= \pf} }
& \binom{\mf}{\mf'} \binom{\pf}{\pf'} 
\la \tauf^{\mf' + \delf_k + \delf_l + \delf_0} \kaf^{\pf'} \ra_0 \,
\la \tauf^{\mf'' + \delf_0} \kaf^{\pf'' + \delf_{a-1}} \ra_0.
\end{align*}
Equivalently, for each $k,l\ge0$ the function $H_0(\tf;\sf)$
satisfies
\begin{align*}
&\dd_a \dd_k \dd_l H_0 = 
(\dd_k \dd_l \dd_0 H_0) (\dd_{a-1} \dd_0 H_0) 
\quad \text{when $a\ge 1$}, \\
& d_1 \dd_k \dd_l H_0 = 
(\dd_k \dd_l \dd_0 H_0) ((\EE-1) \dd_0 H_0), \\
& d_a \dd_k \dd_l H_0 = 
(\dd_k \dd_l \dd_0 H_0) (d_{a-1} \dd_0 H_0)
\quad \text{when $a\ge 2$}.
\end{align*}
This system together with $H_0(t_0,\zef;\zef)= \frac{t_0^3}{6}$ uniquely
determines $H_0$. \qed
\end{thm}

\begin{proof}
The recursion relations are a direct consequence of \eqref{zero:psi},
\eqref{zero:kappa}, and the restriction properties of the $\psi$ and
$\kappa$ classes described above. 

In order to derive the differential equations from the recursion
relations one notices that the increment of $m_a$ or $p_a$ by one in a
recursion relation corresponds to taking the partial derivative with
respect to $t_a$ or $s_a$.

The operator $\EE$ appears because the second recursion relation when
$a=1$ produces $\omega_0$. The corresponding moduli space is
$\M_{0,J\sqcup *}$. As $|J\sqcup *| = |\!|\mf''|\!|+1$, it follows
that $\omega_0 = |J|-1= |\!|\mf''|\!|-1$, and we use that the
multiplication by $m_i$ can be expressed by $t_i\dd_i$.
\end{proof}

\begin{rem}
Setting $\sf=\zef$ in the first equation one gets differential
equations satisfied by $F_0$ (cf. \cite{W}).
\end{rem}

\begin{rem}
Setting $k=l=0$, $t_0=x$, and $t_1=t_2=\ldots=0$ in the second
equation one gets differential equations satisfied by $H(x,\zef;\sf)$
whose third derivative with respect to $x$ is $K_0 (x;\sf)$. These
equations are a simple consequence of the results in
\cite[Sec.~1]{KMZ}.
\end{rem}

Now we turn to genus $1$. We use the explicit presentations
\eqref{one:psi} and \eqref{one:kappa} and take into account the
automorphism groups of the graphs to obtain the following 

\begin{thm}
\label{H1}
For each $\mf \in \SC_0$, $\pf \in \SC_1$, and $a\ge 1$ one has 
\begin{align*}
\la \tauf^{\mf + \delf_a} \kaf^\pf \ra_1 = 
\frac{1}{24} & \la \tauf^{\mf + 2\delf_0 + \delf_{a-1}} 
	 \kaf^\pf \ra_0 \\
+ \sum_{\substack{\mf'+\mf''=\mf\\ \pf' + \pf''= \pf} }
& \binom{\mf}{\mf'} \binom{\pf}{\pf'} 
\la \tauf^{\mf' + \delf_0} \kaf^{\pf'} \ra_1 
\la \tauf^{\mf'' + \delf_0 + \delf_{a-1}} \kaf^{\pf''} \ra_0, \\
\la \tauf^{\mf} \kaf^{\pf + \delf_a} \ra_1 = 
\frac{1}{24} & \la \tauf^{\mf + 2\delf_0} 
	 \kaf^{\pf + \delf_{a-1}} \ra_0 \\
+ \sum_{\substack{\mf'+\mf''=\mf\\ \pf' + \pf''= \pf} }
& \binom{\mf}{\mf'} \binom{\pf}{\pf'} 
\la \tauf^{\mf' + \delf_0} \kaf^{\pf'} \ra_1 
\la \tauf^{\mf'' + \delf_0} \kaf^{\pf'' + \delf_{a-1}} \ra_0.
\end{align*}
Equivalently, the functions $H_1(\tf;\sf)$ and $H_0(\tf;\sf)$ satisfy
\begin{align*}
& \dd_a H_1 = \frac{1}{24} \dd_{a-1} \dd_0 \dd_0 H_0 +
(\dd_0 H_1) (\dd_{a-1} \dd_0 H_0) 
\quad \text{when $a\ge 1$}, \\
& d_1 H_1 = \frac{1}{24} \EE \dd_0 \dd_0 H_0 +
(\dd_0 H_1) ( (\EE-1) \dd_0 H_0), \\
& d_a H_1 = \frac{1}{24} d_{a-1} \dd_0 \dd_0 H_0 +
(\dd_0 H_1) (d_{a-1} \dd_0 H_0) 
\quad \text{when $a\ge 2$}. 
\end{align*}
This system together with the system and the initial conditions from
Thm.~\ref{H0} uniquely determines the pair $H_0, H_1$. \qed
\end{thm}

\begin{rem}
The first recursion relation when $\pf=\zef$ was obtained by Witten in
\cite{W}.
\end{rem}

The system of differential equations above can be solved explicitly
for $H_1$ in terms of $H_0$ to derive \eqref{our:rel}. 

\begin{crl}
The functions $H_1$ and $H_0$ are related by 
\begin{equation*}
H_1 = \frac{1}{24} \log \dd_0^3 H_0.
\end{equation*}
\end{crl}

\begin{proof}
Because of uniqueness it is enough to check that $\frac{1}{24}
\log \dd_0^3 H_0$ satisfies the differential equation in
Thm.~\ref{H1}. This is a straight forward calculation which makes use
of the differential equations from Thm.~\ref{H0}. 
\end{proof}

\begin{rem}
Setting $\sf=\zef$ we recover a result from \cite[Sec.~2.2]{Di}.
Setting $t_0=x, t_1=t_2=\ldots=0$ we get $K_1 = \frac{1}{24} \log
K_0$.
\end{rem}

%	derivation of the puncture/dilaton equation
% Last Modified 05/30/97 4pm by TK
\section{Puncture and Dilaton Equations}
\label{punc}

In this section we introduce an approach which does not use explicit
presentations of $\psi$ and $\kappa$ classes in terms of graphs.
Instead we introduce the analogues of the puncture and dilaton
equations. These equations generalize the classical puncture and
dilaton equations obtained by Witten \cite{W,W2}. (We shall explain
these equations below.) This will allow us to write differential
equations for $H$, and then, using these differential equations, prove
that $H_0$ and $H_1$ satisfy \eqref{our:rel}.

\subsection{Recursion Relations}

In Sec.~\ref{kmz} we introduce the notation incorporating the
intersection numbers of both of the $\psi$ and $\kappa$ classes. Now
we shall to prove certain recursion relations for these numbers. 

\begin{lm}
\label{punc:dil}
The following recursion relations are satisfied:
\begin{equation}
\label{puncture}
\la \tauf^{\mf+\delf_0} \kaf^\pf \ra = 
\sum_{i=1}^\infty m_i
\la \tauf^{\mf+\delf_{i-1}-\delf_i} \kaf^\pf \ra + 
\sum_{\substack{\jf=\zef\\ |\jf|>0}}^\pf \binom{\pf}{\jf}
\la \tauf^\mf \kaf^{\pf-\jf+\delf_{|\jf|-1}} \ra, 
\end{equation}
and for each $a\ge 1$
\begin{equation}
\label{dilaton}
\la \tauf^{\mf+\delf_a} \kaf^\pf \ra = 
\sum_{\jf=\zef}^\pf \binom{\pf}{\jf}
\la \tauf^\mf \kaf^{\pf-\jf+\delf_{|\jf|+a-1}} \ra.
\end{equation}
\end{lm}

\begin{proof}
We continue to use the notation from \ref{basic}, i.e., $\pi$ is the
universal curve over $\M_{g,n}$, $(\psio_i,\ \kao_i)$ and $(\psi_i,\
\kappa_i)$ are classes upstairs and downstairs respectively,
$\sigma_i$ is the $i^{\text{th}}$ canonical section of $\pi$, and
$D_i$ is its image.

It was shown in \cite{W2} that $\psio_i^a = \pi^* \psi_i^a + \pi^*
\psi_i^{a-1} D_i$ and in \cite{AC} that $\kao_i = \pi^* \kappa_i +
\psio_{n+1}^i$. Note also that $\psio_i D_i = 0$, $\psio_{n+1} D_i =
0$ for $i=1,\ldots,n$, and $D_i D_j = 0$ when $i \neq j$. Using this
one derives that
\begin{align*}
\pi_* (\psio_1^{d_1} \ldots & \psio_n^{d_n} 
\kao_1^{p_1} \kao_2^{p_2} \ldots) \\
&= \pi_* \big(\ (\pi^* \psi_1^{d_1} + \pi^* \psi_1^{d_1-1} D_1) \ldots
         (\pi^* \psi_n^{d_n} + \pi^* \psi_n^{d_n-1} D_n) \\
& \qquad\qquad   \times (\pi^* \kappa_1 + \psio_{n+1})^{p_1}
         (\pi^* \kappa_2 + \psio_{n+1}^2)^{p_2} \ldots \big) \\
&= \sum_{i:d_i \neq 0} 
   \psi_1^{d_1} \ldots \psi_i^{d_i-1} \ldots \psi_n^{d_n} 
   \kappa_1^{p_1} \kappa_2^{p_2} \ldots \\
& \qquad\qquad + \sum_{\substack{\jf=\zef\\ |\jf|>0}}^{\pf}
   \binom{\pf}{\jf} \psi_1^{d_1} \ldots \psi_n^{d_n} 
   \kappa_{|\jf|-1} \kappa_1^{p_1-j_1} \kappa_2^{p_2-j_2} \ldots \\
\intertext{Similarly one can show that}
\pi_* (\psio_1^{d_1} \ldots & \psio_n^{d_n}  \psio_{n+1}^a
\kao_1^{p_1} \kao_2^{p_2} \ldots) \\
& = \sum_{\jf=\zef}^{\pf}
   \binom{\pf}{\jf} \psi_1^{d_1} \ldots \psi_n^{d_n} 
   \kappa_{|\jf|+a-1} \kappa_1^{p_1-j_1} \kappa_2^{p_2-j_2} \ldots
\end{align*}
Recall that $\kappa_0 = 2g-2+n$. 

One can further integrate the push forward formulas above to obtain
the statement of the lemma.
\end{proof}

\begin{rem}
Recursion relations \eqref{puncture} and \eqref{dilaton} do not mix
intersection numbers in different genera.
\end{rem}

\begin{rem}
If $\pf = \zef$, then the second sum in the first relation vanishes,
and we obtain the classical puncture equation. If $\pf=\zef$ and $a=1$
in the second relation, then we get the classical dilaton equation.
Note that both classical equations involve only $\psi$ classes.
\end{rem}

\begin{rem}
This is clear recursion relations \eqref{puncture} and
\eqref{dilaton} allow to eliminate $\tauf$ from the intersection
number, i.e., to express all mixed intersection numbers through the
intersection numbers on $\M_{0,3}$, $\M_{1,1}$, and the intersection
numbers of the $\kappa$ classes on $\M_{g,0}$, $g\ge 2$.

In \cite[Cor.~2.3]{KMZ} the authors obtained an explicit expression
for the intersection numbers of the $\kappa$ classes through the
intersection numbers of the $\psi$ classes. This should be related to
\eqref{dilaton}, but we do not know how to derive their formula from
it.
\end{rem}

\subsection{Differential Operators}

Now we derive differential equations for $H$ using recursions
\eqref{puncture} and \eqref{dilaton}. Recall that $\dd_i$, $d_i$
denote the partial derivatives with respect to $t_i$, $s_i$.

\begin{thm}
The function $\exp (H(\tf;\sf))$ is annihilated by the following
differential operators:
\begin{align*}
-\dd_0 \ + & \  
\sum_{\jf:\, |\jf|\ge 2} \frac{\sf^\jf}{\jf!} d_{|\jf|-1} +
\sum_{i=0}^\infty t_i \dd_{i-1} \\
+ & \ s_1 ( \sum_{i=0}^\infty \frac{2i+1}{3} t_i \dd_i +
\sum_{i=1}^\infty \frac{2}{3} i\, s_i d_i) +
\frac{1}{2} t_0^2 \delta_{g,0} + \frac{1}{24} s_1 \delta_{g,1},\\
-\dd_1 \ + & \ 
\sum_{\jf:\, |\jf|\ge 1} \frac{\sf^\jf}{\jf!} d_{|\jf|} +
(\sum_{i=0}^\infty \frac{2i+1}{3} t_i \dd_i +
\sum_{i=1}^\infty \frac{2}{3} i\, s_i d_i) + 
\frac{1}{24} \delta_{g,1},\\
-\dd_a \ + & \  
\sum_{\jf} \frac{\sf^\jf}{\jf!} d_{|\jf|+a-1} 
\quad \text{when $a\ge 2$.} 
\end{align*}
\end{thm}

\begin{rem}
The differential operators above do not mix genus, and therefore they
annihilate each $\exp (H_g (\tf;\sf) )$ separately. 
\end{rem}

\begin{rem}
The function $\exp (F(\tf))$ is annihilated by differential operators
$L_i$, $i\ge -1$, which, after a rescaling of variables, satisfy the
Virasoro relations \cite{Di}. The first two differential operators in
the statement of the theorem are analogues of $L_{-1}$ and $L_0$
respectively, which encode the puncture and dilaton equations,
respectively.
\end{rem}

\begin{proof}
The differential operators above are the direct translation of the
recursion relations \eqref{puncture} and \eqref{dilaton}. The addition
of $\delf_i$ to $\mf$ or $\pf$ translates into taking the
corresponding partial derivative. The subtraction of $\delf_i$ from
$\mf$, and multiplying the term by $m_i$ translates into the
multiplication by $t_i$. One should also change the summation index to
obtain the second summand of each differential operator.

The terms in parentheses in the first two equations come from the
value of $\kappa_0= 2g-2+n$. We use \eqref{charge} in order to express
this number in terms of differential operators. Finally, the last terms
in the first two equations correspond to the initial conditions $\la
\tau_0^3 \ra_0=1$, $\la \kappa_1 \ra_1 = \frac{1}{24}$, and $\la
\tau_1 \ra_1= \frac{1}{24}$.
\end{proof}

The theorem above leads to another proof of Cor.~\ref{our:rel}. First
one notes that $H_0$ and $H_1$ are uniquely determined by the
differential operators above. Therefore it suffices to check that
$\frac{1}{24} \dd_0^3 H_0$ satisfies the genus $1$ equations. This is
a direct calculation.

%	explicit formulas
% Last Modified 05/30/97 4pm by TK
\section{Explicit Expressions}
\label{explicit}

In this section we write the closed form expressions for the
intersection numbers in genus $1$ as sums of the multinomial
coefficients.  

\begin{nota}
If $\bb= (b_1, \ldots, b_k)$ is a vector with integer entries we
denote by $\lf\bb\rf$ the multinomial coefficient $\frac{(b_1+ \cdots
+b_k)!} {b_1! \ldots b_k!}$, and we set it to zero if at least one
entry is negative. Recall also that $|\!| \bb |\!|$ denotes the sum
$b_1+ \ldots +b_k$. 
\end{nota}

In genus $0$ the intersection numbers of the $\psi$ classes are very
simple: $\la \tau_{b_1} \ldots \tau_{b_k} \ra_0 = \lf \bb \rf$. In
order to state our result in genus $1$ we define for each $k\ge 1$ a
function $f_k: \nz_{\ge 1}^k \to \nz_{\ge 1}$ by
\[
f_k(\bb) = f_k(b_1, \ldots, b_k) := 
\la \tau_0^{|\!| \bb |\!| - k} \tau_{b_1} \ldots \tau_{b_k} \ra_1.
\]
Clearly each $f_k$ is invariant under the permutations of its
arguments. 

\begin{prop}
For each $k\ge1$
\begin{equation}
\label{expl:psi}
f_k (\bb) = \frac{1}{24} \lf \bb \rf - \frac{1}{24}
\sum_{\substack{\epsf \in \{ 0,1 \}^k\\ |\!| \epsf |\!|\ge 2}}
(|\!| \epsf |\!| - 2)!\, \lf \bb - \epsf \rf.
\end{equation}
\end{prop}

\begin{proof}
The intersection numbers of the $\psi$ classes in genus $1$ are
determined by the classical puncture equation and the classical
dilaton equation. (See the second remark after Lemma \ref{punc:dil}.)
Reformulated in terms of the collection $\{ f_k \}$ these equations
say that this collection is uniquely determined by the following
properties:
\begin{itemize}
\item $f_1 (b_1) \equiv \frac{1}{24}$,
\item $f_k$ is invariant under the permutation of the arguments,
\item $f_k (\bb) = \sum_{i=1}^k f_k (\bb - \delf_i)$ when $b_i \ge 2$ for all $i$,
\item $f_k (b_1, \ldots, b_{k-1}, 1) = (b_1+ \cdots +b_{k-1}) f_{k-1}
(b_1, \ldots, b_{k-1})$.
\end{itemize}
The first three properties are obviously satisfied by the expression
given in the statement of the proposition. A direct computation verifies
the last property. 
\end{proof}

An explicit expression for the intersection numbers of the $\kappa$
classes in genus $1$ can be obtained by substitution of
\eqref{expl:psi} into Cor.~2.3 from \cite{KMZ}. This expression is
quite complicated. We do not know how to simplify this expression, and
therefore we do not present the resulting formula for the $\kappa$
classes here. 

%	asymptotics
% Last Changed  5/30/97 4pm by TK
\section{Asymptotic Formulas for Volumes of $\M_{1,n}$}
\label{asymp}

In this section, we derive an asymptotic formula for the
Weil-Petersson volumes of $\M_{1,n}$ in the limit that $n$ becomes
very large extending the proof of a similar result for genus zero in
\cite{KMZ}. We do so by using analytic properties of the generating
function $K(x;\sf)$ in the case where all $s_i=0$ for all $i\geq 2$.
In some sense, these results are complementary to those of Penner
\cite{Pe} who obtains similar formulas for the case where the genus
becomes very large.

The class of the Weil-Petersson symplectic form on $\M_{g,n}$ is precisely
$\frac{1}{2\pi^2} \kappa_{(g,n),1}$. The symplectic volume of $\M_{g,n}$ 
is called the \emph{Weil-Petersson volume of} $\M_{g,n}$. For this reason,
the intersection numbers associated to the $\kappa$ classes are sometimes
called \emph{higher Weil-Petersson volumes}. To avoid unnecessary factors, we
shall work instead with the quantity

\begin{df}
$w_{g,n} := \int_{\M_{g,n}}\kappa_1^{3g-3+n}$.
\end{df}

\begin{thm}[\cite{KMZ}]
The genus zero Weil-Petersson volumes satisfy the asymptotic relation as
$n\,\to\,\infty$ 
\[
w_{0,n+3}\,\thicksim\,
\frac{\gamma_0\, 2^\frac{3}{2}}{C\sqrt{\pi}} \frac{2^{2n} n^{2n +
\frac{1}{2}}}{C^n e^{2n}},
\]
where $\gamma_0\,\approx\,2.40482555777\ldots $ is the smallest zero of
the Bessel function $J_0$ and $C = - 2\gamma_0 J_0'(\gamma_0)\,\approx\,
2.496918339\ldots$.
\end{thm}

They proved this by noticing that the function $H_0''(x,\zef;\sf)$ is
invertible, and when all $s_i=0$ except for $s_1$ the inverse function
satisfies Bessel's equation, after a change of variables. Combining their
results with ours for genus one, we obtain the following.

\begin{thm}
The genus one Weil-Petersson volumes satisfy the asymptotic relation as
$n\,\to\,\infty$ 
\[
w_{1,n}\,\thicksim\,\frac{\pi}{24}\frac{(2n)^{2n}}{C^n e^{2n}},
\]
where $C$ is same constant as in the above. 
\end{thm}

\begin{proof}
One uses the asymptotic formulas for the genus zero case and our
result that the generating functions are related by \eqref{our:rel} .
\end{proof}

\begin{rem}
The theorem above supports a conjecture of Itzykson regarding the
existence of such an asymptotic formula for all genera with a constant
$C$ independent of the genus (cf.~\cite[p.~765]{KMZ}). 
\end{rem}

%	cohomological field theories
% Last Changed by TK 5/30/97 4:30pm
%
\section{The Moduli Space of Cohomological Field Theories}
\label{cft}

The moduli space of normalized, rank one cohomological field theories of
genus zero was described in Kontsevich, Manin, and Zagier \cite{KMZ}. The
generating function associated to the $\kappa$ classes endows this moduli
space with coordinates which behave nicely with respect to taking tensor
products of cohomological field theories, a notion introduced in \cite{KMK}
(see also \cite{Ka}).

In this section, we introduce the notion of a restricted, normalized,
cohomological field theory in genus one and describe the moduli space
of rank one theories of this kind. Such \cfts\ turn out to be almost
completely determined by their genus zero part using the relations
between the boundary strata of $\M_{1,n}$ recently obtained by Getzler
\cite{G}. The analogous set of coordinates are constructed for this
moduli space but to do so, we must introduce the $\lambda$ classes, as
well.

\subsection{Cohomological Field Theories}

Consider $\G_{g,n}$, the set of stable graphs of genus $g$ and $n$
tails labeled with the set $[n]$ . Each $\G_{g,n}$ is acted upon by
the permutation group $S_n$ which permutes the labels on the tails.
There are composition maps
\[ \G_{g_1,n_1} \times \G_{g_2,n_2}\,\to\,\G_{g_1+g_2,n_1+n_2-2}\]
taking $(\Gamma,\Gamma') \, \mapsto \, \Gamma \circ_{(i_1,i_2)}
\Gamma'$ for all $i_1$ in $[n_1]$ and $i_2$ in $[n_2]$ given by
grafting the tail $i_1$ of $\Gamma$ with tail $i_2$ of $\Gamma'$ and
then relabeling the remaining tails with elements of the set
$[n_1+n_2-2]$ by inserting orders. There are another set of
composition maps $\G_{g,n}\,\to\,\G_{g+1,n-2}$ taking
$\Gamma\,\mapsto\,\tr_{(i_1,i_2)}\Gamma$ for all distinct pairs $i_1$
and $i_2$ in $[n]$ in which the tails $i_1$ and $i_2$ of $\Gamma$ are
grafted together. These composition maps are equivariant with respect
to the action of the permutation groups. Let $\GG{g,n}$ be the vector
space over $\nc$ with a basis $\G_{g,n}$ then the compositions and
permutations group actions can be extended $\nc$-linearly. Let $\GG{}$
denote the direct sum of $\GG{g,n}$ for all stable pairs $g,n$.

The collection $\{\,\GG{g,n}\,\}$ (or, for that matter,
$\{\,\G_{g,n}\,\}$) together with the composition maps and actions of
the permutation groups described above forms an example of a
\emph{modular operad}, a notion due to by Getzler and Kapranov
\cite{GK}. By restricting to just the genus zero subcollection
$\{\,\GG{0,n}\,\}$ and forgetting about the composition maps $\tr$, we
obtain an example of a \emph{cyclic operad} \cite{GK2}.

Similarly, the homology groups $H_\bullet(\M_{g,n})$ are endowed with an
action of $S_n$ which relabels the punctures on the stable curve and
there are composition maps
\[ \circ_{i_1,i_2}\,:\,H_{p_1}(\M_{g_1,n_1})\,\otimes\,
H_{p_2}(\M_{g_2,n_2})\,\to\,H_{p_1+p_2}(\M_{g_1+g_2,n_1+n_2-2}) \] for
all $i_1$ in $[n_1]$, $i_2$ in $[n_2]$ and $\tr_{(i_1,i_2)}\,:\,
H_{p}(\M_{g,n})\,\to\,H_{p}(\M_{g+1,n-2})$ for all distinct $i_1$ and
$i_2$ in $[n]$, both of which are induced from the inclusion of
strata. These composition maps are equivariant under the action of the
permutation groups.

The natural maps $\alpha_{g,n}\,:\,\GG{g,n}\,\to\,
H_\bullet(\M_{g,n})$ mapping $\Gamma\,\mapsto\,[\M(\Gamma)]$ where
$\M(\Gamma) := \prod_{v\in V(\Gamma)}\,\M_{g(v),n(v)}$ preserves the
above structures and gives rise to the sequence of morphisms

\begin{equation}
\label{eq:boundary}
0\,\longrightarrow\,\left< \R_{g,n}\right> \,\longrightarrow\,\GG{g,n}\,
\overset{\alpha_{g,n}}{\longrightarrow} \, H_\bullet(\M_{g,n})
\end{equation}
where the kernel of $\alpha_{g,n}$ is denoted by $\left< \R_{g,n}\right>$,
the ideal in $\GG{}$ generated by some space of relations $\R_{g,n}$.

\begin{df}
The modular operad $\HH\, :=\,\{\,\HH_{g,n} \,\}$ is the collection of
\[
\HH_{g,n}\,:= \, \frac{\GG{g,n}}{\left< \R_{g,n}\right>}
\]
\end{df}

The canonical diagonal maps $\M_{g,n}\,\to\,\M_{g,n}\times \M_{g,n}$ induce
maps $\,H_\bullet(\M_{g,n})\,\to\,
H_\bullet(\M_{g,n})\otimes H_\bullet(\M_{g,n})$ making $H_\bullet(\M_{g,n})$
into a Hopf modular operad in the natural way \cite{GK}. This endows
$\HH_{g,n}$ with the structure of a Hopf modular operad, as well.

In the case of $g=0$, the results of \cite{KM} and \cite{Ke} implies that
$\alpha_{0,n}$ is surjective and $\HH_{0,n}$ is  isomorphic to
$H_\bullet(\M_{0,n})$. Furthermore, the relations $\R_{0,n}$ are those due to
Keel \cite{Ke} which come from a lift of the basic codimension one relations
$\R_{0,4}$ on $\M_{0,4}$ via the canonical forgetful map
$\M_{0,n}\,\to\,\M_{0,4}$.

In the case of $g=1$, $\alpha_{1,n}$ is known not to be surjective since
 $\M_{1,n}$ has odd dimensional homology classes.  However, Getzler
\cite{G} has shown that the space of relations $\R_{1,n}$, in addition to
those coming from Keel's relations, contains the lifts of two other
relations. The first is the lift of the basic codimension one relation 
on $\M_{1,2}$ which contains no genus one vertices -- this may be regarded
as the image of Keel's relations under the self-sewing morphism
$\tr_{(3,4)}\,:\,\GG{0,4}\,\to\,\GG{1,2}$. The second relation, which
contains genus one vertices, is between codimension two strata and is of the
form 

\begin{equation}
12 \delta_{2,2} - 4\delta_{2,3} - 2 \delta_{2,4} + 6\delta_{3,4} +
\delta_{0,3} + \delta_{0,4} - 2\delta_\beta = 0
\label{eq:getzler}
\end{equation}
where each term is an $S_4$-invariant combination of graphs of a given
topological type and each graph $\Gamma$ represents the homology class
$[\M_\Gamma]$. (See \cite{G} for details.) Getzler also states \cite{G} that
he has shown \cite{G2} that $\alpha_{1,n}$ maps surjectively onto the even
dimensional homology of $\M_{1,n}$ and that the relations mentioned above
do in fact generate all of $\R_{1,n}$.

A cohomological field theory  is essentially a representation, in the sense
of operads, of $H_\bullet(\M_{g,n})$. In order to define such an object, we
need to define the appropriate notion of the endomorphisms of a vector space
is in this context. Let $V$ be a vector space over $\nc$ with a symmetric,
nondegenerate bilinear form $h$ of degree zero. Let $\End{V}_{g,n}:=T^n V$ be
the $n^{\rm th}$ tensor power of $V$ for all nonnegative integers $g,n$ such
that $2g-2+n>0$ where $T^0 V$ is understood to be $\nc$. $S_n$ acts upon 
$\End{V}_{g,n}$ by permuting the tensor factors and the composition maps
$\End{V}_{g_1,n_1} 
\otimes \End{V}_{g_2,n_2}\,\to\,\End{V}_{g_1+g_2,n_1+n_2-2}$ taking
$(\mu,\mu') \,  \mapsto \, \mu \circ_{(i_1,i_2)} \mu'$ for all $i_1$ in
$[n_1]$ and $i_2$ in $[n_2]$ given by applying the inverse of $h$ to the
the corresponding tensor factors of $\mu$ and $\mu'$, and inserting the
remaining factors in the usual way. Similarly, the composition
$\End{V}_{g,n}\,\to\,\End{V}_{g+1,n-2}$ taking
$\mu\,\mapsto\,\tr_{(i_1,i_2)}\mu$  for all distinct pairs $i_1$ and  $i_2$
in $[n]$ corresponds to applying the inverse of $h$ to the appropriate pair
of tensor factors of $\mu$.

\begin{df}[Cohomological Field Theory]
A \emph{(complete) cohomological field theory (\cft) of rank} $r$, $(V,h)$,
is a morphism of modular operads
$\mu_{g,n}\,:\,H_\bullet(\M_{g,n})\,\to\,\End{V}_{g,n}$ where $(V,h)$ is an
$r$-dimensional vector space with an invariant, symmetric bilinear form
$h$. A \emph{\cft\ of genus $g$ } are maps $\mu_{g',n}\,:\,
H_\bullet(\M_{g',n})\,\to\,\End{V}_{g',n}$ which are defined only for 
$g'\leq g$ which satisfy all the axioms of a \cft\  in which no higher genus
maps appear. A \emph{restricted \cft\ } is a morphism
$\mu_{g,n}\,:\,\HH_{g,n}\,\to\,\End{V}_{g,n}$.
\end{df}

A \cft\ can also be described described dually in terms of maps
$\End{V}_{g,n}\, \to\, H^\bullet(\M_{g,n})$. 

Notice that a restricted \cft\ of genus zero is the same as a \cft\ of genus
zero since $H_\bullet(\M_{0,n})\,=\, \HH_{0,n}$.

\begin{rem}
In the language of \cite{W2,W}, a \emph{topological gravity (coupled to
topological matter)} is a \cft\ and the morphisms $\mu_{g,n}$ are
the \emph{correlation functions} of the theory. The genus zero \cft\ is said
to be \emph{tree level} while a genus one \cft\ is said to be \emph{one
loop}.
\end{rem}

\begin{rem}
The natural Hopf structure on $H_\bullet(\M_{g,n})$ endows the category of
\cfts\  with a tensor product as is usual in representation theory.
\end{rem}

An restricted \cft\ is completely determined by a generating function called
its potential. If the \cft\ is not restricted then one can still define the
notion of a potential (essentially since the modular operad
$H_\bullet(\M_{g,n})$ is the quotient of some free modular operad) but
we will not need to work in such generality.

\begin{df}
The \emph{potential} $\Phi = \sum_{g=0}^\infty \Phi_g$ of a restricted \cft\
$\mu\,:\,\HH\,\to\,\End{V}$ of rank $r$ is defined by choosing a basis
$\{\,e_1,\ldots,e_r\,\}$ for $V$ where
$I_{g,n}(e_{a_1},e_{a_2},\ldots,e_{a_n})$ is the number obtained by using $h$
to pair $\mu_{g,n}([\M_{g,n}])$ with $e_{a_1}\otimes
e_{a_2}\otimes\ldots\otimes e_{a_n}$ and 
\[
\Phi_g(\mathbf{x}) := \sum_{n=0}^\infty
\,I_{g,n}(e_{a_1},e_{a_2},\ldots,e_{a_n}) \,\frac{x^{a_1} x^{a_2}\ldots
x^{a_n}}{n!}. 
\]
(where the summation convention has been used) which is regarded as an
element in $\nc[[x^1,\ldots,x^r]]$.
\end{df}

\begin{thm}
\label{thm:wdvvg}
A element $\Phi_0$ in $\nc[[x^1,\ldots,x^r]]$ is the potential of a
rank $r$, genus zero \cft\ $(V,h)$ if and only if \cite{KM,Ma} it
satisfies the WDVV equation
\[
(\partial_{a} \partial_{b} \partial_{e}\Phi_0)\, h^{ef}\, (\partial_{f}
\partial_{c} \partial_{d}\Phi_0)\, = \,(-1)^{|x_a|(|x_b| + |x_c|)}\,
(\partial_{b} \partial_{c}
\partial_{e} \Phi_0)\, h^{ef} \,(\partial_{f} \partial_{a}
\partial_{d}\Phi_0),
\]
where $h_{a,b} := h(e_a,e_b)$, $h^{ab}$ is in inverse matrix to $h_{ab}$,
$\partial_a$ is derivative with respect to $x^a$, and the summation
convention has been used.

If $(\Phi_0,\Phi_1)$ is the potential associated to a restricted, rank $r$
\cft\ of genus one then $\Phi_0$ must satisfy the WDVV equation and
$(\Phi_0,\Phi_1)$ must satisfy Getzler's equation from proposition
$(3.14)$ in \cite{G}. \qed
\end{thm} 

The WDVV equation can be read off from the basic codimension one relation on
$\M_{0,4}$. Similarly, Getzler's equation can be seen from his relation
(equation \ref{eq:getzler}). The second statement will become an if and only
if after the proof in \cite{G2} appears.

\subsection{Rank One Cohomological Field Theories}

Let $(V,h)$ be a rank one \cft\ with a fixed unit vector $e$. The
morphisms $\HH_{g,n}\,\to\,\End{V}_{g,n}$ are completely determined by
the collection of numbers $\{\,I_{g,n}\,\}$ where $I_{g,n} :=
\mu_{g,n}([\M_{g,n}])(\underbrace{e,e,\ldots,e}_n)$ which must satisfy
relations between themselves reflecting the way that the boundary
strata in $\M_{g,n}$ fit together. The potential in this case is
\[
\Phi_g = \sum_{n=0}^\infty I_{g,n} \frac{x^n}{n!},
\]
where $I_{g,n}$ is defined to vanish for pairs $(g,n)$ which are not
stable.

We will see that tautological classes on the moduli space of curves
give rise to complete rank one \cfts. In order to describe the moduli
space of restricted, rank one \cfts\ of genus one, we need to introduce a
combination of the $\lambda$ classes which behave nicely with respect
to restriction.

\begin{df}
For all stable pairs, $(g,n)$, let $\Lambda_{g,n}$ be an element in
$H^\bullet(\M_{g,n})[\mathbf{s},\mathbf{u}]$ (where $\mathbf{s} =
(s_1,s_2,\ldots)$ and $\mathbf{u} = (u_1,u_2,u_3,\ldots\,)$) then
let 
\[
\Lambda_{g,n} := \exp(\,\sum_{i=1}^\infty (\,s_i\,\kappa_{(g,n),i}\, +\,
u_{i}\,\gamma_{(g,n),i}\,)\,)
\]
where $\gamma_{(g,n),i} := \mathrm{ch}_{2i-1}( \pi_*\,\omega_{g,n})$. Here
$\mathrm{ch}_i$ is the $i^{\mathrm{th}}$ Chern character, and
$\pi_*\,\omega_{g,n}$ is the pushforward of the relative dualizing
sheaf.  (Notice that $\mathrm{ch}_{2i}( \pi_*\,\omega_{g,n})$ vanishes for all
$i$ \cite{Ma,Fa}.) 
\end{df}

The classes $\gamma_{(g,n),i}$ are polynomials in the $\lambda$ classes.  In
particular, $\gamma_1 = \lambda_1$. 

\begin{thm}
The collection $\Lambda := \{\,\Lambda_{g,n}\,\}$ gives rise to a
complete, rank one \cft\ for all values of $\textbf{u}$ and
$\textbf{s}$ by integrating the cohomology classes $\Lambda_{g,n}$ over the
homology classes on $\M_{g,n}$. Furthermore, the tensor product of the \cft\
associated to parameter values $(\mathbf{s_1},\mathbf{u_1})$ and
$(\mathbf{s_2},\mathbf{u_2})$  is the \cft\  associated to
$(\mathbf{s_1}+\mathbf{s_2},  \mathbf{u_1}+ \mathbf{u_2})$. Similarly,
\end{thm}

\begin{proof}

In the case where $u_i$ vanishes for all $i$, this was proven in \cite{KMK}
where it was realized that the $\kappa$ classes form a logarithmic \cft\ 
following the work of Arbarello and Cornalba \cite{AC} (see equation
\ref{eq:lcft}).

The proof for the case where all the $s_i$ vanish is as follows.
Consider the bundles $E_{g,n} := \pi_*\, \omega_{g,n}$ on $\M_{g,n}$
from section \ref{tautological}. If $\Gamma$ is a graph of genus $g$
with $n$ tails, then it determines the morphism $\rho_\Gamma:
\prod_{v\in V(\Gamma)} \M_{g(v),n(v)} \to \M_{g(\Gamma), S(\Gamma)}$.
The pull back of $E_{g,n}$ under $\rho_\Gamma$ differs from
$\oplus_{v\in V(\Gamma)} E_{g(v),n(v)}$ by a trivial bundle. It
follows that for each $k\ge 1$ the collection of the Chern characters
$\gamma_k = ch_{2k-1} E_{g,n}$ forms a logarithmic \cft\ (see
\cite{Fa,Fa2}). 

The first part of the theorem follows by combining these two results. The
proof that the coordinates $(\mathbf{s},\mathbf{u})$ are additive with
respect to taking tensor products follows from the definition of coproduct
which is induced from the diagonal map.
\end{proof}

\begin{crl}
The potential of the rank one \cft\ associated to $\Lambda$ (for given values
of $\mathbf{s}$ and $\mathbf{u}$) is precisely the generating 
function $\chi_g$ for the intersection numbers of $\kappa_i$ and
$\gamma_{i}$ classes
\[
\chi_g(x;\mathbf{s},u) := \left< \exp(x\,\tau_0\,+\, 
\sum_{i=1}^\infty\,(\,s_i\,\kappa_i + u_{i}\,\gamma_{i}\,))\right>_g =
\sum_{n=0}^\infty\,I_{g,n}\, \frac{x^n}{n!}\] 
where
\[ I_{g,n} = \sum_{\mathbf{r},\mathbf{m}}\,
\,\frac{\mathbf{s}^\mathbf{m}}{\mathbf{m}!}\,\frac{\mathbf{u}^\mathbf{r}
}{\mathbf{r}!}  
\left<\,\kaf^\mathbf{m}\,\tau_0^n\,\boldsymbol{\gamma}^\mathbf{r}\,\right>_g.
\]
It is understood that $I_{g,n}\,:=\,0$ for unstable pairs $(g,n)$, 
\end{crl}

Notice that the \cft\ arising  $\Lambda$ have the property that $I_{0,3} =
1$ for all values of $\mathbf{u}$ and $\mathbf{s}$. This motivates the
following definition which will play an important role in what 
follows. 

\begin{df}
A rank one, \cft\ of genus $g$ is said to be \emph{invertible} if $I_{0,3}$
is nonzero and \emph{normalized} if $I_{0,3}=1$. 
\end{df}

\subsection{Cohomological Field Theories in Genus Zero and One}

Let us recall the results of Kaufmann, Manin, and Zagier for rank one \cfts\
of genus zero \cite{KMZ}. A rank one \cft\ in genus zero is uniquely
determined by its potential $\Phi_0(x) = \sum_{n=3}^\infty
I_{0,n}\frac{x^n}{n!}$. Furthermore, any function $\Phi_0(x)$  in
$x^3\,\nc[[x]]$ arises from some rank one,  \cft\ of genus zero since the WDVV
equation is trivially satisfied for rank one theories. Therefore, the moduli
space of \cfts\ of genus zero are parameterized by the independent variables
$I_{0,n}$ for $n\geq 3$. What is nontrivial, however, is the behavior of
these potentials under tensor product. In particular, the coordinates
$I_{0,n}$ do not behave nicely under tensor product. However, the generating
function associated to the genus zero $\kappa$ classes
$H_0(x,\mathbf{0};\mathbf{s})$ (which is equal to $\chi_0$ with $u=0$) allows
them to introduce coordinates on the space of normalized, rank one 
\cfts\  which behave nicely under tensor products. 

\begin{thm}[\cite{KMZ}]
\label{thm:kmz}
The moduli space of normalized, rank one \cfts\ in genus zero are
parameterized by $\mathbf{s}$ with potential $\Phi_0(x;\mathbf{s}) =
H_0(x,\mathbf{0};\mathbf{s})$ in $\nc[\mathbf{s}][[x]]$, our generating
function for the intersection numbers of $\kappa$ classes in genus zero.
Furthermore, taking tensor products is additive with respect to the
coordinates $\mathbf{s}$.
\end{thm}

We now treat the case of genus one and discover that the $\kappa$ classes
are not sufficient to describe the entire moduli space of normalized, rank
one \cfts.  We will see that one needs to introduce the first $\lambda$
classes. 

\begin{thm}
If the pair $(\Phi_0,\Phi_1)$ is a potential associated to a restricted, rank
one \cft\ of genus one  then the following equation holds in $\nc[[x]]$
\[
\label{eq:rk1getzler}
-(\Phi_0^{(3)})^2\Phi_1^{(2)} + \Phi_0^{(3)} \Phi_0^{(4)}
\Phi_1^{(1)} - \frac{1}{12} (\Phi_0^{(4)})^2 +
\frac{1}{24}\Phi_0^{(3)}\Phi_0^{(5)} = 0 
\]
where $\Phi_g^{(l)}$ is the $l^{\rm th}$ derivative of $\Phi_g$.
\end{thm}

\begin{proof}
Equation \ref{eq:rk1getzler} above is nothing more than the equation due to
Getzler in theorem \ref{thm:wdvvg} for the case of rank one theories. Our
equation can be seen from equation \ref{eq:getzler} directly by associating
to each graph 
\[
\Gamma\,\mapsto\,
\frac{1}{|\textrm{Aut}(\Gamma)|} \prod_{v\in V(\Gamma)}\,\frac{\partial^{n(v)}
\Phi_{g(v)}}{\partial x^{n(v)}}
\]
and then extending linearly to linear combinations of graphs. One will obtain
$-36$ times the equation above. 
\end{proof}

Unlike the case of genus zero where the WDVV equation is trivially satisfied,
solutions to this equation fall into two classes depending upon whether
$\Phi_0^{(3)}$ is invertible in the ring of formal power series $\nc[[x]]$. 

\begin{thm}
\label{thm:rk1solu}
The pair $(\Phi_0,\Phi_1)$ is a potential associated to an invertible, restricted,
rank one \cft\  of genus one if and only if $\Phi_0(x)$ is of the form
$I_{0,3}\,\frac{x^3}{6} + x^4\,\nc[[x]]$ for $I_{0,3}$ nonzero and
\[
\label{eq:rk1solu}
\Phi_1 = \frac{1}{24}\,\log\Phi_0''' + B \Phi_0''
\]
where $B$ is an arbitrary constant. Therefore, an invertible, restricted \cft\ of
genus one is uniquely determined by arbitrary values of $I_{0,n}$ for all
$n\geq 4$, $I_{0,3}\,\not=\,0$, and $I_{1,1}$.  

If the restricted, rank one, \cft\ is not invertible then $\Phi_0\,=\,0$ and
$\Phi_1$ obeys no constraints. Therefore, the space of such theories is
parameterized by all values of $I_{1,n}$ for all $n\geq 1$.
\end{thm}
\begin{proof}

If the pair $(\Phi_0,\Phi_1)$ is a potential associated to an invertible
restricted, rank one \cft\  then since $\Phi_0'''$ has an inverse in $\nc[[x]]$,
one can solve equation \ref{eq:rk1getzler} explicitly.

The converse is more difficult in the absence of the proof that the lifts of
the relations described above genus $1$ span the entire space of relations
$\R_{1,n}$.  However, we will not need this statement but will explicitly
construct restricted, normalized, rank one \cfts\ in genus one which realize all
solutions to equation \ref{eq:rk1solu} above. This we do in the next
subsection.

Since $I_{1,1} = I_{0,4} + B I_{0,3}$, when $I_{0,3}$ nonzero varying $B$
is the same as varying $I_{1,1}$ and leaving all of the $I_{0,n}$ unchanged.

In the case that the restricted, rank one \cft\ is not invertible then our result
follows from the equation.
\end{proof}

We shall not discuss noninvertible \cfts\ any further in this paper.  From
now on, we shall restrict ourselves to normalized \cfts.

It is worth observing that by incorporating the $\psi$ classes, one can use
the previous result to obtain yet another proof of the formula $H_1 =
\frac{1}{24} \log H_0^{'''}$. It is not clear which of these approaches will
prove most useful in higher genera.

\subsection{Potentials in Genus Zero and One}

In this subsection, we construct potentials for a class of normalized, 
restricted, rank one \cfts\ in genus one explicitly and show that they span the
entire space of solutions to equation \ref{eq:rk1solu} completing the
proof of that theorem. These potentials are generating functions associated
to the $\kappa$ classes and $\lambda_1$. This will give rise to coordinates
which are additive under tensor product in analogy with the case of genus
zero in \cite{KMZ}.

We begin with a useful lemma.

\begin{lm}
The tautological class $\lambda_1$ on $\M_{1,n}$ can be written in terms of
boundary classes as follows:
\begin{equation*}
\lift{1500}{$\lambda_1\,=\, \frac{1}{12}$}
\ \begin{picture}(0,0)%
\includegraphics{./pic/lambda.pstex}%
\end{picture}%
\setlength{\unitlength}{0.00033300in}%
\begingroup\makeatletter\ifx\SetFigFont\undefined
% extract first six characters in \fmtname
\def\x#1#2#3#4#5#6#7\relax{\def\x{#1#2#3#4#5#6}}%
\expandafter\x\fmtname xxxxxx\relax \def\y{splain}%
\ifx\x\y   % LaTeX or SliTeX?
\gdef\SetFigFont#1#2#3{%
  \ifnum #1<17\tiny\else \ifnum #1<20\small\else
  \ifnum #1<24\normalsize\else \ifnum #1<29\large\else
  \ifnum #1<34\Large\else \ifnum #1<41\LARGE\else
     \huge\fi\fi\fi\fi\fi\fi
  \csname #3\endcsname}%
\else
\gdef\SetFigFont#1#2#3{\begingroup
  \count@#1\relax \ifnum 25<\count@\count@25\fi
  \def\x{\endgroup\@setsize\SetFigFont{#2pt}}%
  \expandafter\x
    \csname \romannumeral\the\count@ pt\expandafter\endcsname
    \csname @\romannumeral\the\count@ pt\endcsname
  \csname #3\endcsname}%
\fi
\fi\endgroup
\begin{picture}(2487,1516)(1492,-1874)
\put(3979,-1264){\makebox(0,0)[lb]{\smash{\SetFigFont{8}{9.6}{rm}$[n]$}}}
\put(3751,-586){\makebox(0,0)[lb]{\smash{\SetFigFont{8}{9.6}{rm}$\phantom{1}$}}}
\end{picture}
 
\end{equation*}
\end{lm}
\begin{proof}
The proof follows from the fact that $\lambda_{(1,n),1} =
\pi^*\lambda_{(1,1),1}$ via the forgetful map
$\pi:\,\M_{1,n}\to\,\M_{1,1}$. One uses \eqref{lambda:kappa} to
express $\lambda_{(1,1),1}$ in terms of boundary classes.
\end{proof}

In the sequel, let $\chit_g(\mathbf{s},u)$ be equal to the generating
function $\chi_g(\mathbf{s},\mathbf{u})$ where all values of $u_i$ are set to
zero except for $u\,:=\,u_1$.

\begin{thm}
\label{thm:lambda}
The intersection numbers above satisfy the following:
\[
\chit_0(x;\mathbf{s},u) = H_0(x,\mathbf{0};\mathbf{s})
\]
and
\[
\chit_1(x;\mathbf{s},u) = \frac{u}{24}  H_0''(x,\mathbf{0};\mathbf{s}) +
\frac{1}{24} \log  H_0'''(x,\mathbf{0};\mathbf{s})
\]
where $'$ denotes differentiation with respect to $x$. 
\end{thm}
\begin{proof}

Using that $\lambda_1$ vanishes on $\M_{0,n}$, the presentation of
$\lambda_1$ on $\M_{1,n}$ in terms of boundary strata above,  and  
the fact that the $\kappa$ classes and $\lambda_1$ forms a logarithmic \cft,
we obtain the equations

\[
\left< \kaf^\mathbf{m}\lambda_1^r\tau_0^n\right>_0 = 0.
\]

and

\[
\left< \kaf^\mathbf{m}\lambda_1^r\tau_0^n\right>_1 = 
\begin{cases}
\frac{1}{24}\left<\kaf^\mathbf{m}\tau_0^{n+2}\right>_0 & \text{if
$r=1$,} \\ 0 & \text{if $r\geq 2$.}
\end{cases}
\]

Rewriting these identities in terms of $\chit_g$, using the fact that
$\chit_g(\mathbf{s},x,0) = H_g(\mathbf{s},x)$ and theorem \ref{our:rel}, we
obtain the desired result.
\end{proof}

\begin{proof}{(completion of theorem \ref{thm:rk1solu})}
By setting $B := \frac{u}{24}$ in the previous theorem and $\Phi_g = \chit_g$
for $g=0,1$, we conclude the proof of theorem \ref{thm:rk1solu} since theorem
\ref{thm:kmz} implies that by forgetting $\Phi_1$, we obtain all possible
\cfts\ in genus zero. Furthermore, by varying $u$, one obtains all possible
values of $I_{1,1}$ without changing the values of $I_{0,n}$.
\end{proof}

\begin{rem}
The relations between the intersection numbers obtained in the previous proof
can be encoded in the differential equations
\begin{equation}
\frac{\partial}{\partial u} \chit_0 = 0\qquad\textrm{and}\qquad 
\frac{\partial}{\partial u} \chit_1 = \frac{1}{24} \frac{\partial^2}{\partial x^2}
\chit_0 
\end{equation}
\end{rem}

Putting everything together, we arrive at the following theorem.

\begin{thm}
The moduli space of normalized, restricted, rank one \cfts\ of genus one is
parameterized by coordinates $(\mathbf{s},u)$ via potentials $(\chit_0,\chit_1)$
where $\chit_0(x)$ belongs to $\frac{x^3}{6}+x^4\,\nc[\mathbf{s},u][[x]]$ and
$\chit_1(x)$ belongs to $x\,\nc[\mathbf{s},u][[x]]$ satisfying theorem
\ref{thm:lambda}. The tensor product is additive in the coordinates
$(\mathbf{s},u)$.
\end{thm}

Given two rank one, normalized \cfts\  in genus zero, it is not obvious
how to write down the potential of the tensor product \cft\ explicitly
in terms of the potentials of the tensor factors. In \cite{KMZ}, the authors
show that the operation of tensor product corresponds to multiplication of
the formal Laplace transforms of the two potentials associated to the tensor
factors. Because a rank one, normalized, restricted \cfts\ in genus one is
determined by its genus zero potential and the value of $u$, an explicit
expression for the potential associated to the tensor product of two such
theories follows from the genus zero result of \cite{KMZ} and theorem
\ref{thm:lambda}.

\bibliographystyle{amsplain}
%\bibliography{op}

\begin{thebibliography}{10}

\bibitem{AC}
E.~Arbarello, M.~Cornalba, \emph{Combinatorial and algebro-geometric
cohomology classes on the moduli space of curves,}
J. Algebraic Geom. \textbf{5} (1996), 705--749.

\bibitem{CH}
L.~Caporaso, J.~Harris, \emph{Counting plane curves of any genus},
\texttt{alg-geom/9608025}.

\bibitem{Co}
M.~Cornalba, \emph{On the projectivity of the moduli spaces of
curves,} 
J.~Reine Angew. Math. \textbf{443}, (1993) 11--20. 

\bibitem{DR}
P.~Deligne, M.~D.~ Rapoport,
\emph{Les sch\'emas de modules de courbes elliptiques,}
Proc. Internat. Summer School, Antwerp 1972, 
LNM \textbf{349} (1973), 143--316. 

\bibitem{Di}
R.~Dijkgraaf, \emph{Intersection theory, integrable hierarchies and
topological field theory,}
New Symmetry Principles in Quantum Field Theory, G.~Mack Ed. 
Plenum (1993), 95-158.
  
\bibitem{Du}
B.~Dubrovin, \emph{Geometry of 2D topological field theories,} 
Integrable systems and Quantum Groups, Lecture Notes in Math. {\bf 1620},
Springer, Berlin (1996).
 
\bibitem{E}
T.~Eguchi, K.~Hori, C.~S.~Xiong, \emph{Quantum cohomology and the Virasoro
algebra,} \texttt{hep-th/9703086}.

\bibitem{Fa}
C.~Faber, \emph{A conjectural description of the tautological ring of the
moduli space of curves,}
University of Amsterdam preprint (1996).

\bibitem{Fa2}
\bysame, \emph{Algorithms for computing intersection numbers on moduli
spaces of curves, with an application to the class of the locus of
Jacobians,} 
\texttt{alg-geom/9706006}. 

\bibitem{FP}
W.~Fulton, R.~Pandharipande, \emph{Notes on stable maps and quantum
cohomology,} \texttt{alg-geom/9608011}.

\bibitem{G}
E.~Getzler, \emph{Intersection theory on $\M_{1,4}$ and elliptic
Gromov-Witten invariants,} MPI Preprint No. $96-161$, December 1996, {\tt
alg-geom/9612009}, to appear in J. Amer. Math. Soc.

\bibitem{G2}
\bysame, \emph{Generators and relations for $H_\bullet(\M_{1,n},\nq)$,}
in preparation.

\bibitem{GK} 
E.~Getzler, M.~Kapranov, \emph{Modular operads},
Preprint, Department of Mathematics, MIT, August 1994, {\tt
dg-ga/9408003}, to appear in Compositio Math.

\bibitem{GK2}
\bysame, \emph{Cyclic operads and cyclic homology,}
to appear in ``Geometry, Topology, and Physics for Raoul,'' (ed.
S.~T.~Yau), International Press, Cambridge, MA, 1994.

\bibitem{HL}
R.~M.~Hain, E.~Looijenga, \emph{Mapping class groups and moduli spaces of
curves,} \texttt{alg-geom/9607004},
to appear in Proceedings of the AMS conference on Alg. Geom, Santa
Cruz 1995. 
  
\bibitem{Hi}
N.~Hitchin, \emph{Frobenius manifolds,}
Cambridge University Preprint.

\bibitem{KS}
V.~G.~Kac, A.~Schwartz, \emph{Geometric interpretation of the partition
function of 2D gravity,} 
Phys. Lett. \textbf{257} (1991), 329--334. 

\bibitem{Ka}
R.~Kaufmann, \emph{The intersection form in $H^\bullet(\M_{0,n})$ and the
explicit K\"unneth formula in quantum cohomology},
Internat. Math. Res. Notices (1996), 929--952.

\bibitem{KMZ}
R.~Kaufmann, Yu.~I.~Manin, D.~Zagier, \emph{Higher Weil-Petersson volumes
of moduli spaces of stable n-pointed curves,}
Commun. Math. Phys. \textbf{181} (1996), 763--787.

\bibitem{Ke}
S.~Keel, \emph{Intersection theory of moduli spaces of stable $n$-pointed
curves of genus zero}, 
Trans. AMS \textbf{330} (1992), 545--574.

\bibitem{Ko}
M.~Kontsevich, \emph{Intersection theory on the moduli space of curves
and the matrix Airy function,}
Commun. Math. Phys. \textbf{147} (1992), 1--23.

\bibitem{KM}
M.~Kontsevich, Yu.~I.~Manin, \emph{Gromov-Witten classes, quantum cohomology, and
enumerative geometry},
Commun. Math. Phys. \textbf{164} (1994), 525--562.

\bibitem{KMK}
M.~Kontsevich, Yu.~I.~Manin (with Appendix by R. Kaufmann), \emph{Quantum
cohomology of a product}, 
Invent. Math. \textbf{124} (1996), 313--340.
  
\bibitem{Lo}
E.~Looijenga, \emph{Intersection theory on Deligne-Mumford
compactifications [after Witten and Kontsevich]},
S\'eminaire Bourbaki, Vol. 1992/93, Ast\'erisque \textbf{216} (1993),
187--212.

\bibitem{Ma}
Yu.~I.~Manin, \emph{Frobenius manifolds, quantum cohomology, and
moduli spaces (Chapters I, II, III)},
MPI Preprint No. $96-113$, January 1996.

\bibitem{Ma2}
\bysame, Private communication.

\bibitem{Mat}
M.~Matone, \emph{Nonperturbative model of Liouville gravity,}
J. Geom. Phys. \textbf{21} (1997), 381--398.

\bibitem{Mu}
D.~Mumford, \emph{Towards an enumerative geometry of the moduli space
of curves}, in ``Arithmetic and Geometry,'' (eds. M.~Artin and J.~Tate),
Part II, Progress in Math., Vol. 36,
Birkh\"auser, Basel (1983), 271--328.

\bibitem{Pe} 
R.~C.~Penner, \emph{Weil Petersson volumes},
J. Diff. Geom. \textbf{35} (1992), 559--608.

\bibitem{RT}
Y.~Ruan, G.~Tian, \emph{A mathematical theory of quantum cohomology,} 
J. Diff. Geom. \textbf{42} (1995), 259--367.

\bibitem{W2}
E.~Witten, \emph{On the structure of the topological phase of two-dimensional
gravity}, Nucl. Phys. B \textbf{340} (1990), 281--332.

\bibitem{W}
\bysame, \emph{Two-dimensional gravity and intersection theory on moduli
space}, Surveys in Diff. Geom. \textbf{1} (1991), 243--310.

\bibitem{Z} 
P.~Zograf, \emph{The Weil-Petersson volume of the moduli
spaces of punctured spheres}, in ``Mapping Class Groups and Moduli
Spaces of Riemann Surfaces,'' (eds. R.~M.~Hain
and C.~F.~B\"odigheimer), Contemp. Math. \textbf{150} (1993),
267--372.

\end{thebibliography}
%\end{document}

\providecommand{\bysame}{\leavevmode\hbox to3em{\hrulefill}\thinspace}

\end{document}